\documentclass[reqno]{amsart}[standalone]
\usepackage[notocbasic]{nomencl}
\makenomenclature
\usepackage{booktabs}
\usepackage{units}
\usepackage{amssymb}
\usepackage{amsmath}
\usepackage{mathrsfs}
\usepackage[foot]{amsaddr}
\usepackage{amsfonts}
\usepackage{bm}
\usepackage{amsmath}
\usepackage[square,sort,comma,numbers]{natbib}
\usepackage{graphicx}
\usepackage{xcolor}
\usepackage{caption}
\usepackage{subcaption}
\usepackage{mathbbol}
\usepackage{hyperref}
\usepackage{lipsum}
\usepackage{mathtools}

\newcommand{\expnumber}[2]{{#1}\mathrm{e}{#2}}
\graphicspath{{Images/}}
\usepackage[top=1in, bottom=1.25in, left=1.25in, right=1.25in]{geometry}
\usepackage [autostyle, english = american]{csquotes}
\MakeOuterQuote{"}

\definecolor{vargreen}{rgb}{0.0, 0.5, 0.0}

\usepackage{tikz}
\usetikzlibrary{babel}
\usepackage{etoolbox} %
\usepackage{listofitems} %
\usepackage{pgfmath}
\tikzset{>=latex} %
\colorlet{myred}{red!80!black}
\colorlet{myblue}{blue!80!black}
\colorlet{mygreen}{green!60!black}
\colorlet{mydarkred}{myred!40!black}
\colorlet{mydarkblue}{myblue!40!black}
\colorlet{mydarkgreen}{mygreen!40!black}
\tikzstyle{node}=[very thick,circle,draw=myblue,minimum size=22,inner sep=0.5,outer sep=0.6]
\tikzstyle{connect}=[->,thick,mydarkblue,shorten >=1]
\tikzset{ %
  node 1/.style={node,mydarkgreen,draw=mygreen,fill=mygreen!25},
  node 2/.style={node,mydarkblue,draw=myblue,fill=myblue!20},
  node 3/.style={node,mydarkred,draw=myred,fill=myred!20},
}
\def\nstyle{int(\lay<\Nnodlen?min(2,\lay):3)} %
\begin{document}
\title[]{A physics-based reduced order model for urban air pollution prediction}
\author{Moaad Khamlich$^{1}$}
\author{Giovanni Stabile$^{2,1}$}
\author{Gianluigi Rozza$^{1}$}
\author{L\'{a}szl\'{o} K\"{o}rnyei$^{3}$}
\author{Zolt\'{a}n Horv\'ath$^{3}$}
\address{$^1$ mathLab, Mathematics Area, SISSA, via Bonomea 265, I-34136 Trieste, Italy}
\address{$^2$ Department of Pure and Applied Sciences, Informatics and Mathematics Section, University of Urbino Carlo Bo, Piazza della Repubblica, 13, I-61029, Urbino, Italy}
\address{$^3$ Sz\'echenyi Istv\'an University, Department of Mathematics and Computational Sciences, Gyor-Moson-Sopron, Gyor 9026, Hungary.}
\maketitle
\begin{abstract}
This article presents an innovative approach for developing an efficient reduced-order model to study the dispersion of urban air pollutants. The need for real-time air quality monitoring has become increasingly important, given the rise in pollutant emissions due to urbanization and its adverse effects on human health.
The proposed methodology involves solving the linear advection-diffusion problem, where the solution of the Reynolds-averaged Navier–Stokes equations gives the convective field. At the same time, the source term consists of an empirical time series.

However, the computational requirements of this approach, including microscale spatial resolution, repeated evaluation, and low time scale, necessitate the use of high-performance computing facilities, which can be a bottleneck for real-time monitoring.
To address this challenge, a problem-specific methodology was developed that leverages a data-driven approach based on Proper Orthogonal Decomposition with regression (POD-R) \nomenclature{POD-R}{POD with Regression} coupled with Galerkin projection (POD-G) \nomenclature{POD-G}{POD with Galerkin projection} endorsed with the discrete empirical interpolation method (DEIM). The proposed method employs a feedforward neural network to non-intrusively retrieve the reduced-order convective operator required for online evaluation. The numerical framework was validated on synthetic emissions and real wind measurements.

The results demonstrate that the proposed approach significantly reduces the computational burden of the traditional approach and is suitable for real-time air quality monitoring. Overall, the study advances the field of reduced order modeling and highlights the potential of data-driven approaches in environmental modeling and large-scale simulations.


\end{abstract}

 \printnomenclature
\section{Introduction}
\label{sec:intro}
Urban air pollution is a major global challenge responsible for climate,
ecosystems, and health damage. In particular, in most European Community countries, urban traffic is the most important source of pollutants such as nitrogen oxides, carbon monoxide, particulate matter, and benzene \cite{report-air}. It is believed that in the future, due to continued urbanization and expansion of urban areas, an increasing proportion of the population will be more exposed to high concentrations of traffic pollutants.
Recent research by the World Health Organization (WHO) has shown that 9 out of 10 people breathe polluted air \cite{report-who}. This problem
translates into an economic cost estimated at 5 trillion dollars a year, but
even more important are the consequences on people's health. Research
over the past decades has consistently indicated that air pollution causes
significant damage to the health of exposed populations, and evidence indicates that pollution generated by vehicular traffic is a significant contributor to the adverse effects. The most recent epidemiological studies have shown that the risk from exposure to traffic pollutants is not evenly distributed in a given urban area but increases with decreasing distance from emission sources \cite{zhu2002concentration}. In populations residing in areas of high traffic
volume and density, increased mortality from natural and respiratory causes, increased occurrence of atherosclerosis of coronary and carotid arteries, and increased incidence of ischemic heart disease and bronchial asthma in children have been found
\cite{pope2002lung,brook2004air,hoffmann2007residential,brugge2007near}. These latter medical complications are reflected in human life costs, estimated at 1.8 million deaths globally in 2019 \cite{southerland}.

For these reasons, air quality management has been listed among the United Nations Sustainable Growth Goals, and a European Commission directive has mandated air quality measurement through appropriate monitoring stations and mathematical modeling tools \cite{eu-directive}.

In particular, a framework that combines direct measurements and computational
modeling techniques is a critical analytical tool that can extract various insights
from the collected statistics, for example, a deterministic relationship between
concentrations and emissions and the effectiveness of prevention strategies.

Urban air pollution problems require the study of tracer dispersion in the air.
Indeed, the underlying model is that of atmospheric dispersion
\cite{seinfeld2016atmospheric}, which is a system of
advection-diffusion-reaction partial differential equations (PDEs) coupled
through the chemical production term.   Typically, such problems are addressed
by computational fluid dynamics (CFD) techniques.
In particular, since pollutant dispersion depends on daily weather conditions at the urban scale, CFD models with low time scales, repeated evaluation, and fine mesh discretization must be used. The former requirements translate into
huge memory and computational power, making it essential to use HPC facilities to get
results in reasonable time frames \cite{lazlo21}.

In the present work, we decided to address the criticality represented by
computational cost through the employment of Reduced Order Models (ROMs) \nomenclature{ROM}{Reduced Order Model}
\cite{benner2017model, benner2017, Grepl2007,schilders,vol1,vol2} to achieve
fast converged solutions with limited loss of accuracy.

To simplify the modeling we decided to describe the evolution of the pollutant through the transport equation, where the convective field is given by the solution of the Navier-Stokes equation, while the source term consists of an empirical time series.

We studied two different options for the reduced order model, namely extracting a proper orthogonal decomposition (POD) \nomenclature{POD}{Proper Orthogonal Decomposition} basis onto which the full order empirical source field is projected or using the Discrete Empirical Interpolation Method (DEIM) \nomenclature{DEIM}{Discrete Empirical Interpolation Method} as a hyper-reduction strategy \cite{Chaturantabut2010,Barrault2004,Stabile2020}. Both
these approaches are proven effective, even when the basis for the source term is extracted on a subset of the time series and then used for future state prediction.
We then tackled the parametrized convective field case by changing the
direction and intensity of the inlet velocity. This modeling choice agrees with the aforementioned assumption of coupling the use of our model
with experimental measurements. Here we propose a novel data-driven approach based on a POD-NN \cite{ubbiali} reconstruction of the flux field, which is used to recover in a non-intrusive fashion the reduced-order operators required for the online evaluation.
Our framework is validated on a computational domain modeling the main campus of the University of Bologna, using a mesh with about 40k cells.  We use real inlet conditions for the wind flow around the buildings, based on a one-year long measurement station data with a resolution of one hour. Instead, the source term was obtained synthetically, through a realistic traffic flow modeling, then used to calculate NOx emission.
The developed framework uses tools developed in the integrated urban air pollution dispersion modelling framework developed by Horv\'ath et al. \cite{Horvath2016410} for generating the mesh, the traffic emissions, boundary conditions and the offline simulations.

The work is organized as follows:
\begin{itemize}
    \item Section \ref{sec:related} provides an overview of air pollution modeling and reduction techniques;
    \item Section \ref{sec:problem_formulation} introduces the mathematical problem;
    \item Section \ref{sec:numerical_approximation} provides an in-depth description of the proposed resolution pipeline;
    \item Section \ref{sec:numerical} shows the numerical results for the Bologna campus test case.
\end{itemize}
\section{Related Works}\label{sec:related}
Air pollution models are numerical tools that describe the relationships between atmospheric concentrations, meteorological conditions and emissions. This study can assess the relative impact of the various processes involved, arriving at a deterministic analysis of factors and causes. Air pollution modeling represents a well-established field, the numerical results of which are used in various disciplines: meteorology, engineering, and geography, to name a few.
Primary applications include impact assessment of emission sources, environmental forecasting, inverse modeling, and uncertainty quantification.
For a comprehensive review of the topic, we refer the interested reader to \cite{britter2003flow,sportisse2007review,blocken2011application}.

Pollutant concentrations are usually evaluated using air dispersion models. The problem's difficulty lies in the different scales involved: local, regional, continental, and global.
A first distinction among the modeling techniques employed is that between statistical empirical methods and deterministic methods \cite{zhang2012real}.
However, this distinction is increasingly blurred. In fact, despite being insufficient on their own, measurement stations can supplement and improve the information for deterministic models.
For this reason, a growing interest is directed toward data assimilation techniques, in which an attempt is made to use data from a generic monitoring network to improve predictions and reduce uncertainties related to model input parameters. Some examples of work in this field are \cite{lopez2021urban,hammond2019pbdw,nguyen2021data}.

Like other complex problems, air quality modelling suffers from many uncertainties associated with inputs, a large number of parameters, and numerical discretization. For this reason, much recent work has focused on assessing the sensitivity of model outputs to inputs treated as random variables. Some examples of such works are \cite{carmichael1997sensitivity,he2000application}. In particular, in \cite{sandu2005adjoint}, the study is addressed through the adjoint problem for air quality models.

Another important line of research concerns the more recent use of machine learning methods for pollution prediction.
Typically, machine learning is used as a post-processing tool or directly integrated within the model under consideration \cite{xi2015comprehensive}.
In particular, the recent success of machine differentiation techniques has made it possible to employ neural networks for the regression task in this field as well \cite{boznar1993neural,gardner1998artificial}.

It must be kept in mind that an end-use of air pollution models always requires a drastic reduction in computational costs. Approaches to effect this reduction are among the most varied. For example, early attempts involved simplified physics modeling. In \cite {hanna1971simple}, a method is proposed that estimates the surface concentration from a simple relationship with the production term and wind intensity. Other models of simplified physics involve the use of the singular perturbation theory for splitting the time scales of the problem, as presented in \cite{neophytou2004reduced}.
The computational complexity of this type of study justifies the consideration of test cases with simplified geometries. For example, in \cite{balczo2015flow}, flow and dispersion phenomena are analyzed in an environment represented by a square surrounded by rectangular plan geometries, supplemented by wind tunnel tests and tracer concentration measurements. Their analysis leads to the identification of the main flow structures present in the studied domain.

More recently, order reduction has involved using POD and the reduced basis method. In this field, we mention the works \cite{hammond2019pbdw,hammond2017reduced}, which operate the reduction in a nonintrusive manner, considering a test case representing residential urban pollution. In particular, the results of \cite{hammond2017reduced} are also compared with an intrusive reduction technique based on the Generalized Empirical Interpolation Method (GEIM) \cite{maday2013generalized}.
Our work presents a novelty in the field in that we exploit an efficient hybrid reduction methodology using a nonintrusive convective field reconstruction used within an intrusive concentration field resolution. Our intrusive framework is supplemented using the DEIM as a hyper-reduction strategy, and it is shown to be efficient when trained on a subset of the time series related to emission data and boundary conditions on velocity.

\section{Problem formulation}%
\label{sec:problem_formulation}
The transport-diffusion equation is a linear partial differential equation, which takes the form:
\begin{equation}
\label{eq-transport}
\frac{\partial c}{\partial t} - \nu \Delta c +  \nabla\cdot( \mathbf{u} c )=f;
\end{equation}
where  $c(\mathbf{x}, t):\Omega  \times[0,T) \subset \mathbb{R}^{n} \times \mathbb{R}^{+} \rightarrow
\mathbb{R^{+}}$ (with n=2,3) is the unknown function, which can be thought of as the concentration of a pollutant such as NOx. Specifically, the quantity
$c(\mathbf{x}, t) d V$ represents the mass present at time $t$ in an
infinitesimal neighborhood of the point $\boldsymbol{x}$. Consistently, the mass
of pollutant present in volume $V$ at time $t$ is given by :
\begin{equation}
\int_{V} c(\boldsymbol{x}, t) d V.
\end{equation}

Equation \eqref{eq-transport} is known as the convection-diffusion equation because of the physical interpretation of the various terms in it. The diffusive term  "$- \nu \Delta c$"  represents the rate of
change in concentration due to the difference between the average value of $c$
in a neighborhood of $\boldsymbol{x}$  and the value at $\boldsymbol{x}$
itself. In particular, the diffusive flux is given by \textit{Fick's law}:
\begin{equation}
    \mathbf{q} = - \nu \nabla c.
\end{equation}
The minus sign is because the pollutant travels from high-concentration areas to low-concentration ones. The diffusivity constant $\nu > 0$  is dependent on the particular pollutant under consideration and is
assumed to be constant throughout the domain.

The term "$\nabla\cdot( \mathbf{u} c )$", on the other hand, models the convective transport effect, that is, the transport of the pollutant due to the motion of
the fluid in which it is immersed. In particular, the velocity field
$\mathbf{u}$ is given by the resolution of the Navier-Stokes system, which will be presented in the next section.
Finally, $f$ represents any source or sink term in the equation, which can account for external influences or processes affecting the concentration of the pollutant.

\subsection{The Navier-Stokes equations}%
\label{sub:the_navier_stokes_equations}
Let $\Omega \subset \mathbb{R}^{n},(n=2,3)$ be a bounded domain and consider the following system of equations, known as the incompressible Navier-Stokes system:
\begin{equation}
\label{NavierStokes}
\begin{cases}
\mathbf{u_t}+ \mathbf{\nabla} \cdot (\mathbf{u} \otimes \mathbf{u})- \mathbf{\nabla} \cdot 2 \mu \mathbf{\nabla^s} \mathbf{u}=-\mathbf{\nabla}p &\mbox{ in } \Omega \times [0,T] ,\\
\mathbf{\nabla} \cdot \mathbf{u}=\mathbf{0} &\mbox{ in } \Omega \times [0,T] ,
\end{cases}
\end{equation}
where $\boldsymbol{\nabla}^{\boldsymbol{s}}
\mathbf{u}=\frac{1}{2}(\boldsymbol{\nabla} \mathbf{u} + \boldsymbol{\nabla}
^{T}\mathbf{u})$, and $\mu$ is the so-called kinematic viscosity. The unknowns
are represented by the state of the fluid $X=(\mathbf{u},p)$, namely its
velocity and pressure normalized over a constant density. In this
particular setting, the latter can be interpreted as a Lagrange multiplier needed
to impose the incompressibility condition.

If we introduce a characteristic
velocity $U_{0}$ and a characteristic length $L$ for the problem under
consideration, we define the Reynolds number as $Re=U_0 L/\mu$, which plays a fundamental role in fluid dynamics by identifying different regimes of motion. To characterize turbulent motion, we use the Reynolds Averaged Navier-Stokes (RANS) \nomenclature{RANS}{Reynolds Averaged Navier-Stokes} equations, obtained by time averaging the Navier-Stokes equations.
The RANS equations involve decomposing each characteristic quantity $\phi$ into a time-averaged term $\overline{\phi}$ and a fluctuation $\phi^{\prime}$.
The Reynolds stress tensor is introduced to model the effects of turbulent fluctuations, adding six additional unknowns to the problem and requiring appropriate turbulence models for solution.
For our work, we use the $k$-${\epsilon}$ model \cite{k-epsilon} , which describes the transport of turbulent kinetic energy and its dissipation rate. The RANS approach can provide sufficient information for many practical applications, allowing for significant computational time reduction.

\section{Numerical Approximation}%
\label{sec:numerical_approximation}
In this section, we will present the ingredients needed to tackle numerically the resolution of Equation \eqref{eq-transport} coupled with System \eqref{NavierStokes}. The starting point will be represented by the so-called full order (or high-fidelity) discretization, which aims at providing a high-accuracy approximation to the solution of PDEs. With this regard, we will present only the Finite Volume Method (FVM) \nomenclature{FV}{Finite Volume}, which has been used for our simulation. However, one could employ our framework and change the underline full order discretization strategy. The mitigation of the computational cost required by the presented methodology will be tackled using a reduced order modeling technique.
\subsection{Full order approximation}%
\label{sub:full_order_approximation}
In this section, we will present a very brief introduction to the Finite Volume
Method. This method uses the integral form of conservation laws, and because of this, each term that compares in the formulation has a precise physical
interpretation. Furthermore, FVM is particularly popular for engineering
applications in CFD because of the ease in its numerical implementation, as can
be testified by the wide variety of both commercial and open-source codes based
on this method (our simulations have been carried out with the open-source
library OpenFOAM \cite{openfoam}).
The main characteristics of this method can be summarized as follows:
\begin{itemize}
    \item the domain is subdivided in \textit{control volumes (CV)}, which are
    also called cells, and the conservation laws are enforced on each volume;
    \item each variable is expressed at the center of each cell;
    \item \textit{interpolation formulas} are used to express the value of the
    variables and their gradients: one needs to approximate volume integrals and
    fluxes;
    \item after the discretization, one obtains, as for the finite element case,
    an algebraic equivalent formulation which is constituted by an equation for
    each cell.
\end{itemize}
We do not want to discuss the detail of FVM for a generic problem; this can be found, for instance, in \cite{mangani}. Specifically, we will specialize the discussion for the resolution of the Navier-Stokes system in a steady setting.
The resulting velocity field will later be used to solve the linear transport problem.

As already said, one could recover FVM by integrating the system in a control
volume. However, we will follow an alternative strategy proposed in
\cite{StabileRozza2018}, which results in obtaining the same result by starting
from the weak formulation of the Navier-Stokes equations, namely:
\begin{equation}
    \label{NavierStokesFV}
\begin{aligned}
&\int_{\Omega} \boldsymbol{\nabla} \cdot(\mathbf{u} \otimes \mathbf{u}) \cdot \boldsymbol{v} \ \mathrm{d} \Omega
-\int_{\Omega} \boldsymbol{\nabla} \cdot 2 \mu \boldsymbol{\nabla}^{s} \mathbf{u} \cdot \boldsymbol{v} \ \mathrm{d} \Omega+\int_{\Omega} \boldsymbol{\nabla} p \cdot \boldsymbol{v} \ \mathrm{d} \Omega
+\int_{\Omega} \boldsymbol{\nabla} \cdot \mathbf{u} q \ \mathrm{d} \Omega=0 \quad\\ &\forall(\boldsymbol{v}, q) \in V \times Q;
\end{aligned}
\end{equation}
with:
$$
\begin{array}{l}
    V=H_{0,\Gamma_{0}}^{1}\left(\Omega ; \mathbb{R}^{n}\right)=\left\{\boldsymbol{v} \in H^{1}\left(\Omega ; \mathbb{R}^{n}\right):\left.\boldsymbol{v}\right|_{\Gamma_0}=0\right\}, \quad \text{where} \quad \Gamma_{0}=\{\boldsymbol{x} \in \partial \Omega: \mathbf{u}(\boldsymbol{x})=0\};     \\
Q=L^{2}(\Omega). \end{array}
$$
We then proceed with a discretization of the spaces $V$ and $Q$, based on an underlying discretization of the domain using a tessellation $\mathcal{T}=\left\{\Omega_{e}\right\}_{e=1}^{N_{FV}}$ of non-overlapping polyhedron (\textit{finite volumes}) such that  $\bigcup_{e=1}^{N_{FV}}\Omega_{e} = \Omega$.
The solution is, therefore, sought in the finite-dimensional space:
\begin{equation*}
    V^{h} \times Q^{h}=\left(\operatorname{span}\left\{\mathcal{I}_{k}(\boldsymbol{x})\right\}_{k =1}^{N_{FV}}\right)^{n+1};
\end{equation*}
where $\mathcal{I}_{k}(\boldsymbol{x})$ is the basis function of each finite volume:
\begin{equation}
\mathcal{I}_{k}(\boldsymbol{x})=\left\{\begin{array}{l}
1 \quad \text {if } \boldsymbol{x} \in \Omega_{k}, \\
0 \quad \text {if } \boldsymbol{x} \in \Omega \backslash \Omega_{k}.
\end{array}\right.
\end{equation}

By using the divergence theorem, one can transform the volume integrals of equation \eqref{NavierStokesFV} into surface integrals, thus reaching:
\begin{equation}
    \label{NSvarFV}
\begin{aligned}
&\sum_{e=1}^{N_{F V}} \boldsymbol{v}_{e}\left(\int_{\partial \Omega_{e}} \boldsymbol{n} \cdot(\mathbf{u} \otimes \mathbf{u}) \ \mathrm{d} \Gamma\right.
\left.-\int_{\partial \Omega_{e}} \boldsymbol{n} \cdot 2 \mu \nabla^{s} \mathbf{u} \ \mathrm{d} \Gamma+\int_{\partial \Omega_{e}} \boldsymbol{n} p \ \mathrm{d} \Gamma\right)
+\sum_{c=1}^{N_{F V}} q_{e}\left(\int_{\partial \Omega_{e}} \boldsymbol{n} \cdot \mathbf{u} \ \mathrm{d} \Gamma\right)=0; \\
&\forall (\boldsymbol{v}_{e},q_{e}) \in V^{h} \times Q^{h};
\end{aligned}
\end{equation}
where $\boldsymbol{n}$ is the unit normal vector outgoing from $\Omega_{e}$ .

It is possible to manipulate the previous expression by using quadrature formulas that express the surface integrals as simple summations involving values assumed by the fields on each cell's face (we will use the subscript $f$ for values at the center of the face). This leads to the following expressions:
\begin{itemize}
    \item \textbf{nonlinear convective term}: $\int_{\partial \Omega_{e}} \boldsymbol{n} \cdot(\mathbf{u} \otimes \mathbf{u}) \ \mathrm{d} \Gamma=\sum_{f} \boldsymbol{S}_{\boldsymbol{f}} \cdot \mathbf{u}_{f} \otimes \mathbf{u}_{f}$;
    \item \textbf{diffusive term}: $\int_{\partial \Omega_{e}} \boldsymbol{n} \cdot 2 \mu \boldsymbol{\nabla}^{\boldsymbol{s}} \mathbf{u} \ \mathrm{d} \Gamma=\int_{\partial \Omega_{e}} \boldsymbol{n} \cdot \mu \boldsymbol{\nabla} \mathbf{u} \ \mathrm{d} \Gamma=\mu \sum_{f} \boldsymbol{S}_{\boldsymbol{f}} \cdot(\boldsymbol{\nabla} \mathbf{u})_{f}$;
    \item \textbf{pressure gradient term} : $\int_{\partial \Omega_{e}} \boldsymbol{n} p \ \mathrm{d} \Gamma=\sum_{f} \boldsymbol{S}_{\boldsymbol{f}} p_{f}$;
    \item  \textbf{incompressibility term} : $\int_{\partial \Omega_{e}} \boldsymbol{n} \cdot \mathbf{u} \ \mathrm{d} \Gamma=\sum_{f=1}^{N_{f}} \boldsymbol{S}_{\boldsymbol{f}} \cdot \mathbf{u}_{\boldsymbol{f}}$;
\end{itemize}
being $\boldsymbol{S_{f}}$ the surface area vector.

The interpolation coefficients obtained in the discretization process are utilized for each finite volume to construct an algebraic system of equations. These equations can be reorganized and expressed in matrix form as:
\begin{equation}
\begin{array}{r}
\boldsymbol{C}(\boldsymbol{u})\boldsymbol{u}+\mu \boldsymbol{A}\boldsymbol{u}+\boldsymbol{B p}=0 \ ;\\
    \boldsymbol{P}\boldsymbol{u}=0 \ ;
\end{array}
\end{equation}
where the vectors $\boldsymbol{u}$ and $\boldsymbol{p}$ represent the discrete evaluations of velocity and pressure at the finite volume centers.
Instead, the various matrices discretize specific terms in Equation \eqref{NSvarFV}. In particular, $\boldsymbol{C}$ discretizes the nonlinear term, $\boldsymbol{A}$ represents the diffusive term, $\boldsymbol{B}$ captures the gradient of the pressure, and $\boldsymbol{P}$ accounts for the incompressibility term.
Various strategies can be put in place to solve the former system; in particular, we used the \textit{Semi-Implicit Method for Pressure Linked Equations} (SIMPLE) \cite{Patankar1972ACP}, for which the \textit{k-th} iteration can be summarized as follow:
\begin{enumerate}
\item Solve the discretized momentum equation to compute the intermediate velocity field.
    \item Solve the Poisson equation for the pressure field.
\item Correct the velocities based on the new pressure field.
\item Repeat till convergence.
\end{enumerate}

In this case, the nonlinear term is resolved using the velocity flux at the previous iteration. Clearly, the former segregated procedure is more general than the FVM as it can be applied to other full order discretization procedures, such as the Finite Element Method or the Finite Difference Method.
\subsection{Reduced order model}
\label{sec:rom}
When dealing with a nonlinear PDE, we may be interested in finding the solution
by varying a parameter $\boldsymbol{\mu} \in \mathbb{R}^{p}$ that encodes
certain properties of the problem, such as its geometry, physical properties, or
as we shall see in Section \ref{sec:numerical}, boundary conditions. Formally we are interested in solving a parameterized problem in the form:
\begin{equation}
    \label{parmStrong}
    F(c(t,\boldsymbol{\mu}); \boldsymbol{\mu})=0 \quad \text{with} \quad t \in [0,T], \boldsymbol{\mu} \in \mathcal{P} \subset \mathbb{R}^{p}.
\end{equation}
The nonlinear solution manifold is then given by the set of solutions for each instance of the parameter vector, that is:
\begin{equation}
    \mathcal{S}= \left\{c(t,\boldsymbol{\mu}) \ \mid \ t \in [0,T], \ \boldsymbol{\mu} \in \mathcal{P}\right\}.
\end{equation}
So far, we have seen how to solve the problem for a given value of the parameter
$\boldsymbol{\mu}$, using a high-order approximation, such as the FVM presented
in Section \ref{sub:full_order_approximation}. However, this approach is
computationally inconvenient since the computational cost scales proportionally
to the number of degrees of freedom $N_{FV}$.
In particular, in the context of many-queries applications, when it is necessary to conduct real-time simulations of complex systems, such computational cost may be unsustainable even with HPC facilities.
We point out that this situation is of strong relevance to air quality
monitoring, as previously discussed in Section \ref{sec:intro}. To overcome the
former limitation, over the years, various reduced-order models (ROMs) have been
proposed \cite{benner2017model, benner2017}, including the \textit{Reduced
Basis} (RB) method \nomenclature{RBM}{Reduced Basis Method} \cite{hesthaven2015certified,patera07:book,Quarteroni2016}
chosen for the present work.

Specifically, this method involves processing information from a computationally
expensive phase (\textit{offline}) to construct a suitable linear basis. The
original full order model (FOM) equations are then projected onto this basis to
obtain a reduced model to solve new instances of Problem
\eqref{parmStrong}. Therefore, a critical point in the application of the RB
method is the selection of a proper basis during the offline phase.
In the present work, we use Proper Orthogonal Decomposition (POD) \cite{POD},
which provides a reduced dimensional representation
of a given dataset.

In order to compute the POD basis, we use the method of snapshots as originally
proposed in \cite{Sirovich1987} and implemented in ITHACA-FV\footnote{The source
code of ITHACA-FV can be found in \url{https://github.com/ITHACA-FV/ITHACA-FV}}.

We introduce a time discretization by dividing the interval $[0,T]$ into $N_{T}$
subintervals of equal length $\Delta t=\frac{T}{N_{T}}$, and define $t_{k}=k
\Delta t, \  0\le k \le N_{T}$. We also introduce a discretization of the
parameter space, $\mathcal{P}_{train} \subset \mathcal{P}$, whose cardinality is
$N_{\mu}$. Problem \eqref{parmStrong} is solved for each $\boldsymbol{\mu}_{m}
\in \mathcal{P}_{train}$, obtaining the full order solution \footnote{Here, we will denote with the subscript ${}_{h}$ the solution
to the full order problem, whereas we will use the subscript $_{rb}$ for all the reduced basis quantities.} in the discretized
time interval $\{c_{h}(t_{k},\boldsymbol{\mu}_{m})\}_{k=0}^{N_{T}} \subset
\mathbb{R}^{N_{FV}}$.

For clarity, we combine the indices for temporal and parametric discretizations
into one index:
\begin{equation}
    c_{i}= c_{h}(t_{k}, \boldsymbol{\mu}_{m}),
\quad \text{where} \quad
    \begin{dcases}
i=(m-1)(N_{T}+1) + (k+1), \\
0\le k \le N_{T},\\
1\le m \le N_{\mu}.
    \end{dcases}
\end{equation}
The former offline stage yields $N_{s}=N_{\mu}\times(N_{T}+1)$ full order
samples, also called \textit{snapshots}. The POD modes are constructed in order
to minimize the error between each snapshot and its projection onto the $N_{rb}$-dimensional  POD
space:
\begin{equation}\label{eq:minimize}
  E_{N_{rb}}=\sum_{i = 1}^{N_{s}}
  ||c_i-\sum_{k = 1}^{N_{rb}}a_i^k\varphi_k| | \quad \forall N_{POD} =1,...,N_{s};
\end{equation}
where the coefficients $a_{i}^{k}$ are obtained through projection:
\begin{equation}
\quad a_i^k = (c_i,\varphi_k)_{L^2(\Omega)} \quad \forall i=1,...,N_{s}; \quad \forall k=1,\ldots,N_{POD}.
\end{equation}
The minimization problem in Equation \eqref{eq:minimize} is equivalent to the
solution of the eigenvalue problem:
\begin{equation}
\bm{CQ} = \bm{Q\lambda},
\end{equation}
\begin{equation}
C_{ij}=( c_i,c_j )\textbf{}_{L^2(\Omega)} \quad \forall i,j = 1,...,N_{s},
\end{equation}
where $\bm C \in \mathbb{R}^{N_{s} \times N_{s}}$ is a square matrix, known as the \textit{correlation} matrix. Furthermore, $\bm{\lambda}$ is a diagonal matrix containing the eigenvalues.
Following, the basis $\varphi_i$  is obtained as:
\begin{equation}
\varphi_i=\frac{1}{\sqrt{\lambda_{ii}}} \sum_{j = 1}^{N_{s}} c_jQ_{ij}.
\end{equation}
Finally, the POD bases are constructed, and the resulting space is
\begin{equation}
L_T = [ \varphi_1,...,\varphi_{N_{rb}}] \in \mathbb{R}^{N_{FV}\times N_{rb}},
\end{equation}
where the cardinality $N_{rb}$ must be tailored to the specific problem.

The POD modes are then used to approximate the solution $c(t,\boldsymbol{\mu}) $
for any new value of the parameter with a linear combination:
\begin{equation}\label{eq:tofmu}
c(t,\boldsymbol{\mu}) \approx \sum_{i=1}^{N_{rb}} a_i (\mu,t) \varphi_i (x),
\end{equation}
where $a_i (\mu,t)$ are the parameter-dependent coefficients and $ \varphi_i (x)$
are the parameter-independent basis functions \cite{Quarteroni2016}.

The coefficients of Equation \eqref{eq:tofmu} are then obtained by performing a Galerkin
projection of the FOM residual onto the space spanned by the POD modes in order
to obtain a reduced residual equation:
\begin{equation}\label{eq:rom}
\bm{M_r} \dot{\bm a} - \nu \bm{B_r} \bm a + \bm{C_r} \bm{a} = \bm{f_r}(\bm{\mu}),
\end{equation}
where each term inside Equation \eqref{eq:rom} is obtained by Galerkin projection:
\begin{equation}
\begin{cases}
(\bm{M_r})_{ij} = (\varphi_i,\varphi_j)_{L_2{(\Omega)}}, \\
(\bm{B_r})_{ij} = (\varphi_i,\Delta\varphi_j)_{L_2{(\Omega)}},\\
(\bm{C_r})_{ij} = (\varphi_i,\nabla \cdot (\mathbf{u} \varphi_j))_{L_2{(\Omega)}}, \\
(\bm{f_r}(\bm{\mu}))_{i} = (\varphi_i,f(\mu))_{L_2{(\Omega)}}.\label{eq:source}
\end{cases}
\end{equation}
Implementing the POD-Galerkin (POD-G) method reduces the problem's dimensionality to $N_{rb}$. %
\subsection{Discrete Empirical Interpolation Method}
\label{sec:deim}
The POD-G explained in the previous section, aims to exploit the
essential low-dimensional dynamics in high-dimensional computations
\cite{Benner2021,Brunton2019}.
However, if it is not possible to formulate an affine decomposition for the source term, then computing $\boldsymbol{f}_{r}(\boldsymbol{\mu})$ still relies on the original full order system \cite{Fu2018}.
One can approximate the source term $f(\mu)$ as a linear
combination of hierarchically chosen functions $\chi_i (x)$:
\begin{equation}\label{eq:deim}
f(\mu) \approx \sum_{i=1}^{N_\text{DEIM}} p_i (\mu) \chi_i (x),
\end{equation}
where $ p_i (\mu)$ are parameter-dependent coefficients
\cite{Hesthaven2015}. Once the basis functions $\chi_i (x)$ are calculated, the
parameter-dependent coefficients must be computed. We use the DEIM proposed
in its original formulation \cite{Chaturantabut2010}, where the basis $\chi_i
(x)$ are calculated using the POD, and the coefficients are determined with a
point-wise evaluation of the function $f$ in the locations identified by the
magic points $\bm{q}_i$.  Equation \eqref{eq:deim} is substituted in
Equation \eqref{eq:source} to predetermine a matrix of reduced source terms:
$$\bm{F_r}=\bm{L}^T\bm{U}(\bm{P^{T}}\bm{U})^{-1}\in\mathbb{R}^{N_{rb}\times
N_\text{DEIM}},$$ where $\bm{L} = [\varphi_i,\dots,\varphi_{N_{rb}}]\in \mathbb{R}^{{N_{FV}\times
N_{rb}}}$ is the matrix containing the POD modes for $c$. Furthermore, each term of
the matrix $\bm{U}(\bm{P^{T}}\bm{U})^{-1}\in \mathbb{R}^{N_{FV}\times N_\text{DEIM}}$ is
evaluated using the standard DEIM procedure as reported in
\cite{Chaturantabut2010}. Equation \eqref{eq:rom} can thus be rewritten as:
\begin{equation}\label{eq:rom-deim}
\bm{M_r} \dot{\bm a} - \nu_T \bm{B_r} \bm a + \bm{C_r} \bm{a} = \bm{F_r}\bm{p}(\mu),
\end{equation}
where $\bm{p}(\mu)$ is a vector of coefficients which contains a point-wise
evaluation of the source term in correspondence to the location identified by
the magic points $\bm q_i$: $$(\bm{p}(\mu))_i = f(\mu, \bm{q}_{i}), \quad
i=1,\dots,N_\text{DEIM}.$$
In practice, the previous procedure consists of interpolating the forcing $f(\mu)$ at the magic points $\{q_{i}\}_{i=1}^{N_{\text{DEIM}}}$  and a subsequent projection onto the reduced basis $\{\varphi_{i}\}_{i=1}^{N_{rb}}$.
\subsection{POD with regression}%
\label{sub:pod_with_inteporlation_podi_}

POD with regression (POD-R) is a technique used in the ROM community to build data-driven ROMs, even when the numerical solver that computed the solutions is unavailable. Its characteristic feature lies in replacing the projection phase described in Section \ref{sec:rom} with a regression.
This allows its application also on data coming from experimental measurements, such as \textit{Particle Image Velocimetry} (PIV) images \cite{semeraro2012}.

POD-R consists of first computing the POD basis $\{
\Psi_{i} \}_{i=1} ^{N_{\phi}}$ from a set of high-fidelity
solutions $\{\phi_i\}_{i=1}^{N_\text{train}}$, and then in building the ROMs by regressing the
solution manifold using the POD basis:
\begin{equation}
    \label{exp-rom}
    \phi_{h}(\boldsymbol{\mu})\approx \phi_{rb}(\boldsymbol{\mu})=\boldsymbol{\phi}_{rb}^{\top}(\boldsymbol{\mu}) \boldsymbol{\Psi}(\boldsymbol{x})=\sum_{i=1}^{N_{\phi}} \phi_{rb,i}(\boldsymbol{\mu}) \Psi_{i}(\boldsymbol{x}).
\end{equation}
Assume we have a training database $ \mathcal{D}=
\{((\boldsymbol{\mu}_{j},\boldsymbol{\phi}_h(\boldsymbol{\mu}_{j}))\}_{j=1}^{N_\text{train}}$
consisting of $N_\text{train}$ pairs of parameter values and their respective
full order solution.
The database $\mathcal{D}$ is generally the same as that employed for constructing the POD basis. Each FOM solution can be projected onto the reduced basis to obtain the respective ROM coefficients. This can be done in matrix form as:
\begin{equation}
    \boldsymbol{\phi}_{rb}= \mathbf{V}^{\top} \boldsymbol{\phi}_{h};
\end{equation}
where  $\mathbf{V} \in \mathbb{R}^{N_{FV} \times N_{\phi}}$ is
the matrix whose columns are the POD modes.
We, therefore, also have a database for the reduced coefficients
$ \mathcal{D}_{r}=\{(\boldsymbol{\mu}_{j},\boldsymbol{\phi}_{rb}(\boldsymbol{\mu}_{j}))\}_{j=1}^{N_\text{train}}$ which will be used to construct a regression $\pi: \mathcal{P} \mapsto \mathbb{R}^{N_{\phi}}$ approximating the function $\mathscr{F}: \mathcal{P} \rightarrow \mathbb{R}^{N_{\phi}}: \boldsymbol{\mu} \mapsto \boldsymbol{\phi}_{rb}$.
The full order solution will then be recovered through the following coordinate change:
\begin{equation}
    \label{reconstruction}
    \boldsymbol{\phi}_{h}\left(\boldsymbol{\mu} \right) =\mathbf{V} \pi\left(\boldsymbol{\mu} \right) .
\end{equation}
Recently, this technique has been very successful in the ROM community, finding expression in several applications that differ in the regression model employed, among which the most popular are Radial basis function (RBF) regression \cite{xiao2015}, Gaussian process regression (GPR) \cite{demo-ortali} and artificial neural networks (ANNs) \nomenclature{ANN}{Artificial Neural Network} \cite{xiao2015}. In the present work, we employed precisely a regression based on ANNs, which will be described in Appendix \ref{sub:artificial_neural_networks}.
\subsection{POD-G+NN}%
\label{sub:podnn-g}
This section presents a novel method for a ROM based on POD-G and POD-NN and suitable for solving Equation \eqref{eq-transport}.
The problem inherent in the nonlinearity of the source term $f(x)$ is addressed in the current project using the DEIM introduced in Section \ref{sec:deim}. However, we are also interested in introducing a parametric dependence of the convective field with respect to the boundary conditions. Such a dependence could again be addressed through a hyper-reduction technique, which would entail developing an additional intrusive ROM model to perform the projection of the RANs equations. Such a model would then prove to be the bottleneck for solving the transport equation.
In the present work, therefore, we decided to use the POD-NN approach, first introduced in \cite{ubbiali}, which is a particular POD-R technique exploiting ANNs to recover the reduced coefficients.
The problem concerns the convective matrix $\boldsymbol{C}_{r}$ present in Equation \eqref{eq:rom}, whose components we recall to be:
\begin{equation}
(\boldsymbol{C})_{ij} = (\varphi_i,\nabla \cdot (\mathbf{u} \varphi_j))_{L_2{(\Omega)}}.
\end{equation}
First, a practical premise related to FV discretization is needed. In fact, the discrete FV approximation calculates the previous integral from the contributions of the various cells, which are transferred to the surface:
\begin{equation}
\label{eq:cr-fv}
(\boldsymbol{C})_{ij} \approx \sum_{e=1}^{N_{FV}} \sum_{f=1} ^{N_{f}} \phi_{f} (\varphi_i)_{f}
(\varphi_j)_{f} S_{f};
\end{equation}
where we have introduced the flux field $\phi= \mathbf{u} \cdot \boldsymbol{n}$, which is defined on $\cup _{e=1}^{N_{FV}}\partial \Omega _{e}$.
We are interested in the case where the velocity field $\mathbf{u}$, and hence the flux $\phi$ is itself a function of the parameter vector $\boldsymbol{\mu}$.
We consider a POD basis $\{\Psi_{k}\}_{k=1}^{N_{\phi}}$ for the field $\phi$, which can be used for approximations related to new instances of the parameter by employing the POD-R approach:
\begin{equation}
\label{eq:phi-expansion}
\phi (\boldsymbol{\mu}) \approx \sum_{k=1}^{N_{\phi}} \pi_{NN,k}\left(\boldsymbol{\mu} \right)\Psi_{k} \ .
\end{equation}
In particular, the notation $\pi_{NN}$ is intended to emphasize the use of an ANN for the regression task (Appendix \ref{sub:artificial_neural_networks}).
Thus, we can substitute the expansion of Equation \eqref{eq:phi-expansion} within Equation \eqref{eq:cr-fv}. Recall that the matrix $\boldsymbol{C}_{r}$ is used in the online phase; consequently, an efficient strategy respecting the \textit{offline-online} division paradigm  corresponds to the introduction of a third-order tensor $\mathbf{\Gamma} \in \mathbb{R}^{N_{rb} \times N_{rb} \times N_{\phi}}$ defined as follows:
\begin{equation}
    \left( \mathbf{\Gamma} \right)_{ijk}= \sum_{e=1}^{N_{FV}} \sum_{f=1} ^{N_{f}} \left( \Psi_{k} \right)_{f} (\psi_i)_{f}
(\psi_j)_{f} S_{f}.
\end{equation}
Consequently, $\boldsymbol{C_{r}}$  can be computed by performing a tensor contraction of $ \mathbf{\Gamma}$:
\begin{equation}
    \label{eq:tensor-contraction}
    (\boldsymbol{C_{r}} ( \boldsymbol{\mu} ))_{ij}= \sum_{k=1}^{N_{\phi}} \pi_{NN,k} \left( \mathbf{\Gamma} \right)_{ijk}.
\end{equation}
This matrix can then be used to solve new parametric instances of Problem \eqref{eq:rom}.
\section{Numerical results}%
\label{sec:numerical}
In this section, we will present the results of our simulations for the transport-diffusion problem. In particular, we offer a study applied to an air pollution scenario in a simplified urban geometry, namely the main campus of the University of Bologna. This kind of investigation is of fundamental applicative importance since geometric and modeling simplification still retains the most crucial information regarding dominant features in complex urban environments \cite{blocken2011application}.
We will also present the results for the reduction strategy implemented through the RB method, whose accuracy has been tested on never seen parametric instances of the solution we aim at approximating.
First, we will start with the case where the convective field is nonparametric and fixed as the source term varies. In that case, we propose two different strategies for treating the nonlinearity of the source term: extracting a POD basis onto which the full-order empirical source field is projected or using the DEIM as a hyper-reduction strategy.

We will then generalize the ROM by treating the case in which the convective field is parameterized on the boundary conditions. We will address the latter case through the novel POD-G+NN technique presented in Section \ref{sub:podnn-g}.

\subsection{Preliminaries}

\begin{figure}[hbt!]
    \centering
    \includegraphics[height=7cm]{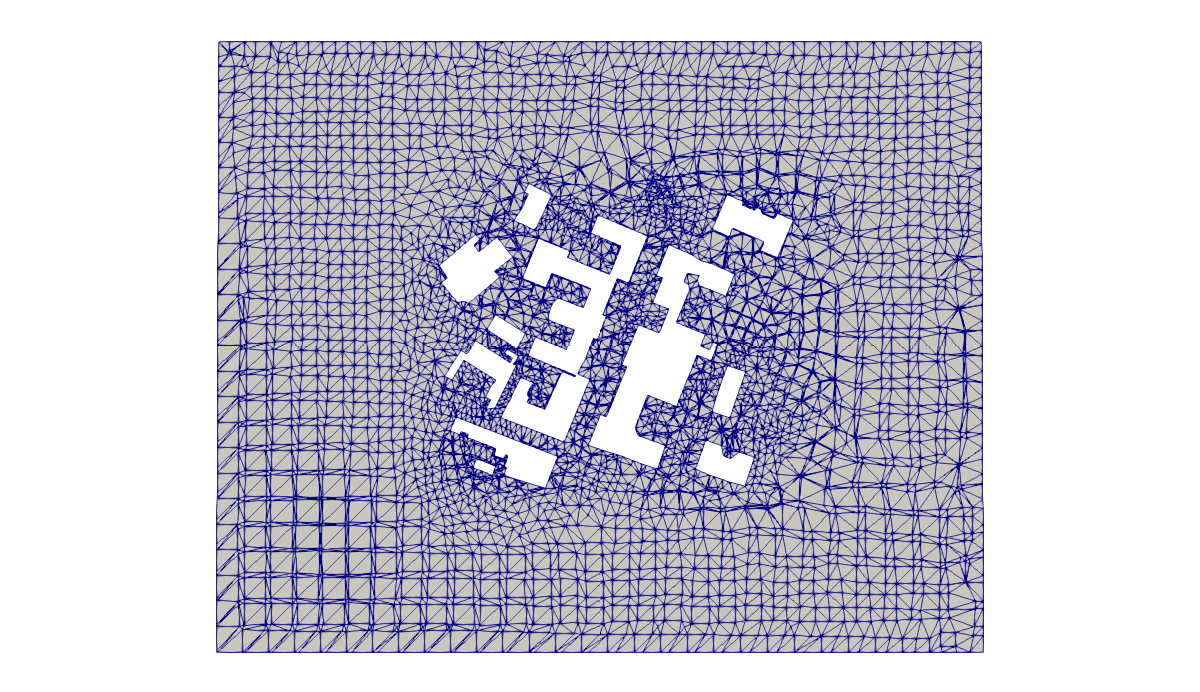}
    \caption{Slice of the computational grid at $z=13\unit{m}$.}
    \label{fig:mesh}
\end{figure}

Let us start by noting that we have used \textit{ SI } unit system for all the equations.
Consequently, the various quantities have the following units of measurement:

\begin {table} [h]
\centering
\begin {tabular} {c|c}
\toprule
\textbf {Physical quantity} & \textbf {Unit of measure}  \\
\midrule
lengths &  \unit{m}\\
times & \unit{s} \\
$c$  & \unitfrac{kg}{m$^3$}  \\
$\mathbf{u}$  & \unitfrac{m}{s}\\
$f$  & \unitfrac{kg}{m$^{3}$s}\\
$ \nu$ & \unitfrac{m$^2$}{s}   \\
$ \mu$ & \unitfrac{m$^2$}{s}   \\
\bottomrule
\end {tabular}
\caption{Unit of measures of the main physical quantities involved in the transport-diffusion problem.}
\end{table}
In the following, we will avoid indicating the units of measurement of the quantities we will be working with to simplify the notation.
As already anticipated in Section \ref{sub:full_order_approximation}, as high-fidelity approximation, we adopted a FV discretization on a fixed 3-dimensional (3-D) tetrahedral unstructured grid with $N_{FV}=37847$ cells. Figure \ref{fig:mesh} shows a slice of the 3-D mesh, where the buildings of the main campus of the University of Bologna are located.
The numerical resolution was addressed using the OpenFOAM library \cite{openfoam}, with the supplementary use of ITHACA-FV \cite{ithaca}, which implements several reduced-order modeling techniques for parameterized problems.
\begin{figure}[hbt!]
    \centering
    \begin{subfigure}{0.45\textwidth}
        \centering
        \includegraphics[trim=100 0 0 100,clip,width=\textwidth]{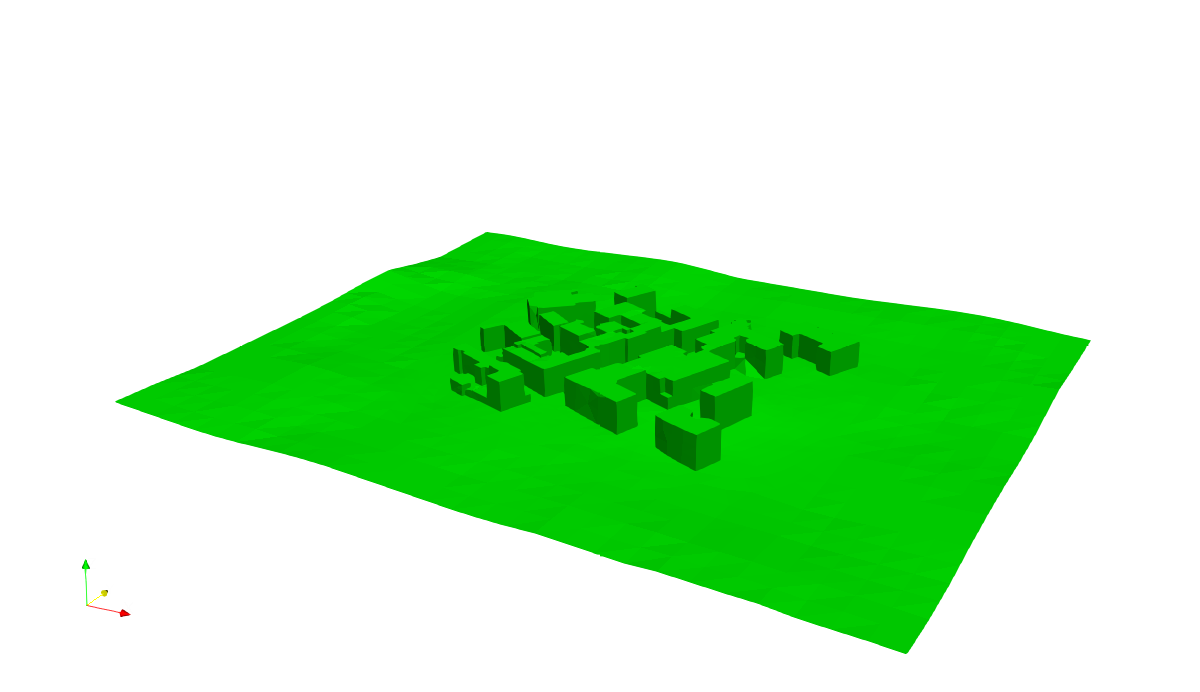}
        \caption{Ground boundary}
    \end{subfigure}
    \begin{subfigure}{0.45\textwidth}
        \centering
        \includegraphics[trim=100 0 0 100,clip,width=\textwidth]{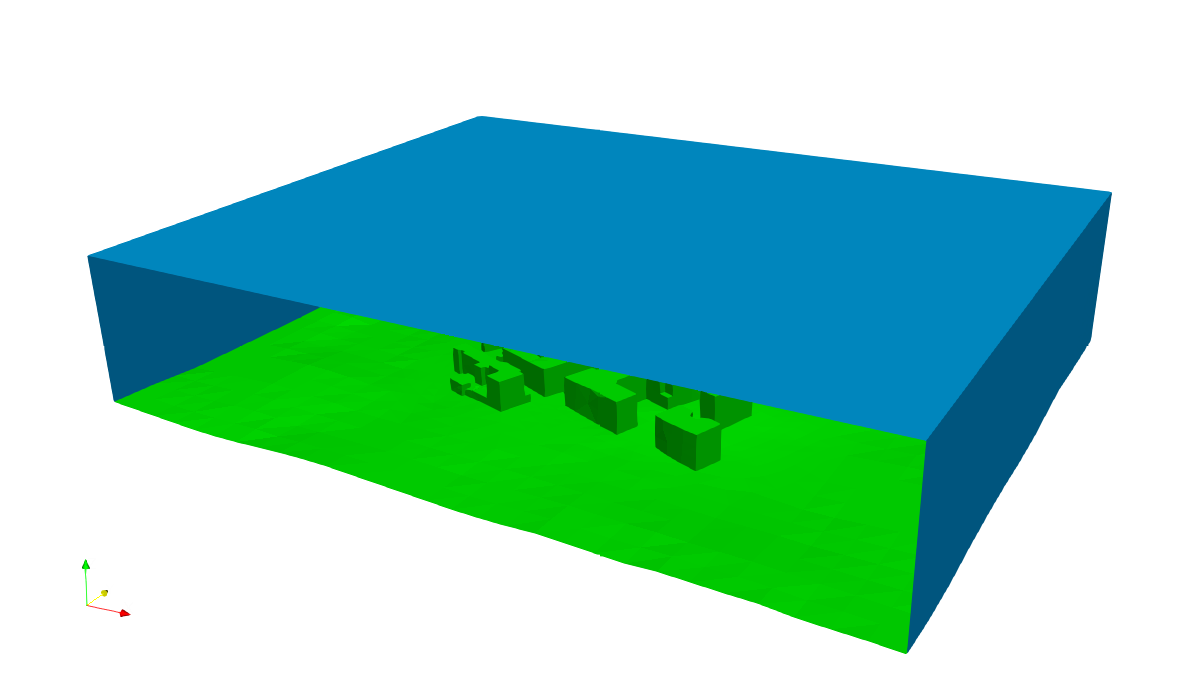}
        \caption{Front and Back boundary}
    \end{subfigure}
    \begin{subfigure}{0.45\textwidth}
        \centering
        \includegraphics[trim=100 0 0 100,clip,width=\textwidth]{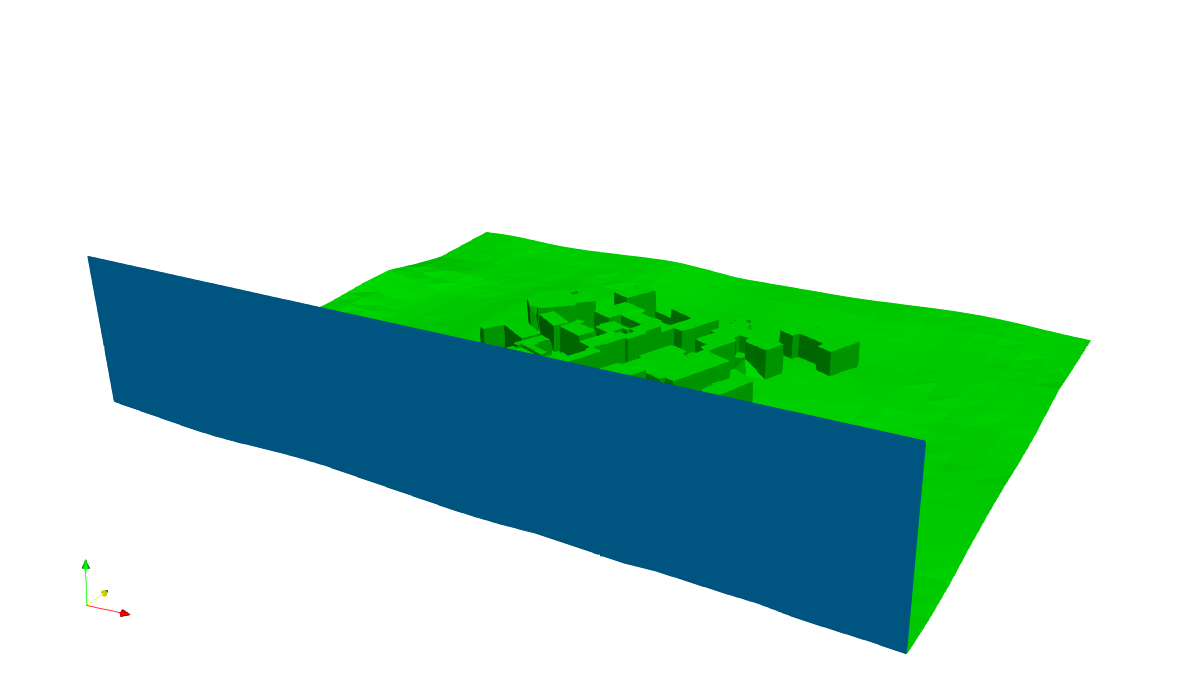}
        \caption{Inlet boundary}
    \end{subfigure}
    \begin{subfigure}{0.45\textwidth}
        \centering
        \includegraphics[trim=100 0 0 100,clip,width=\textwidth]{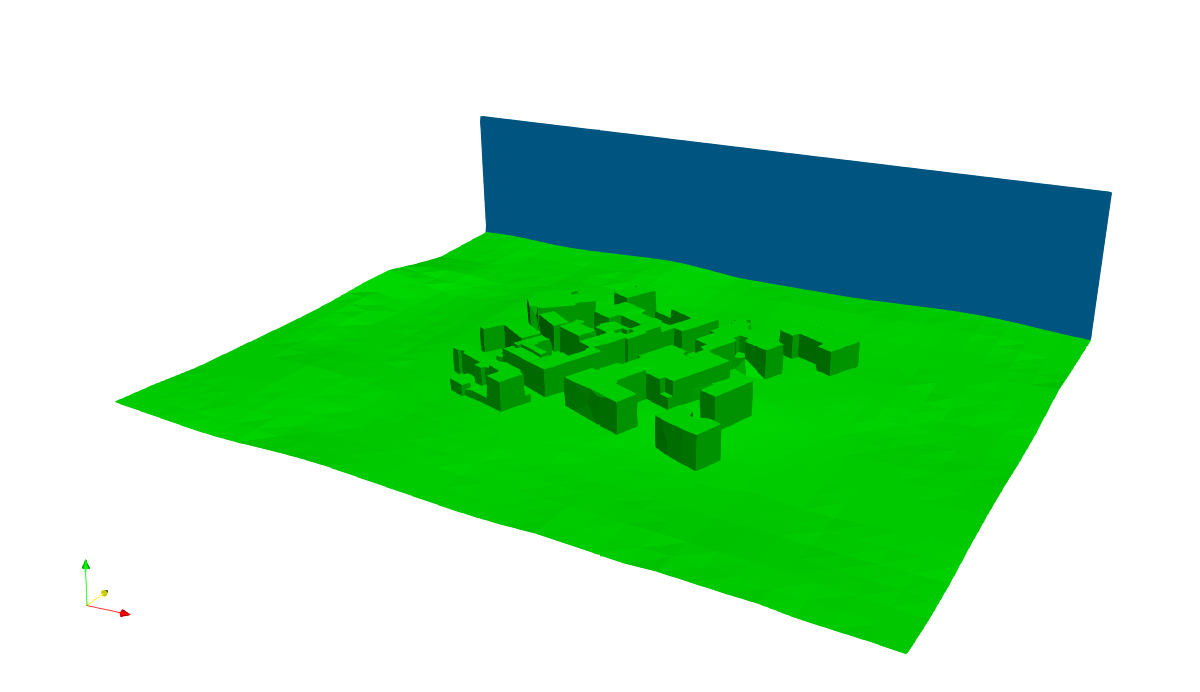}
        \caption{Outlet boundary}
    \end{subfigure}
    \caption{Boundaries of the computational domain.}
    \label{fig:boundaries}
\end{figure}

\subsection{Constant convective field}%
\begin{figure}[hbt!]
    \centering
    \includegraphics[height=7.5cm]{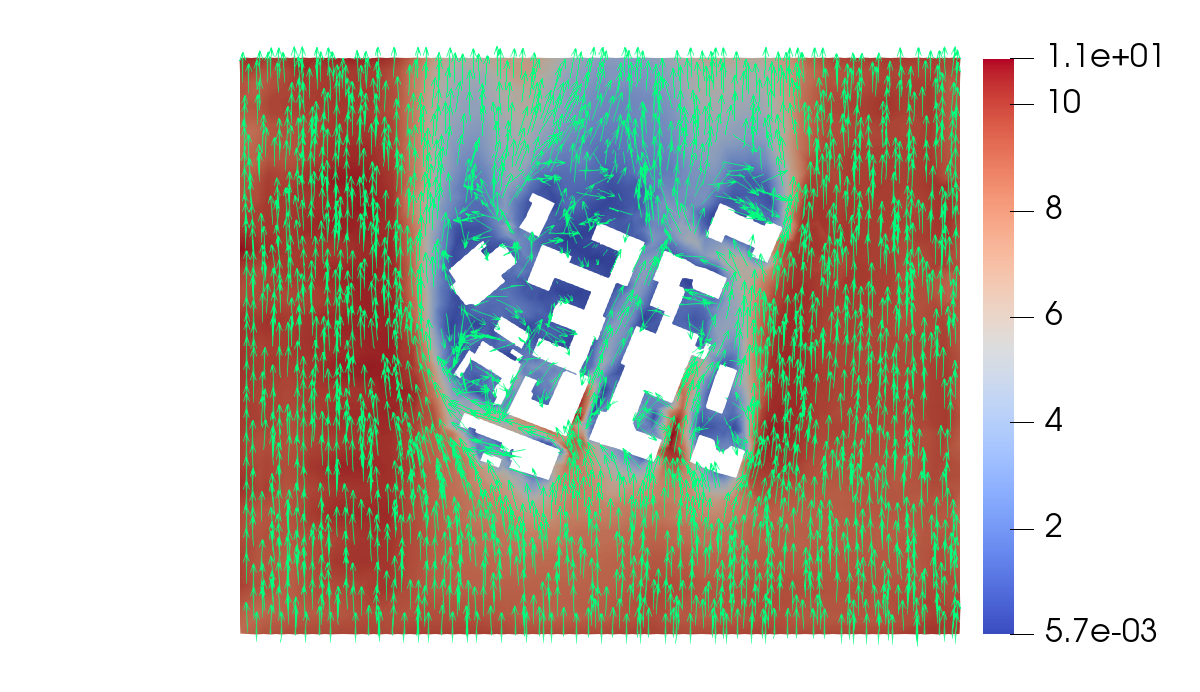}
    \caption{Steady-state velocity field at $z=13m$.}
    \label{fig:velField}
\end{figure}

\begin{figure}[hbt!]
    \centering
    \includegraphics[height=7cm]{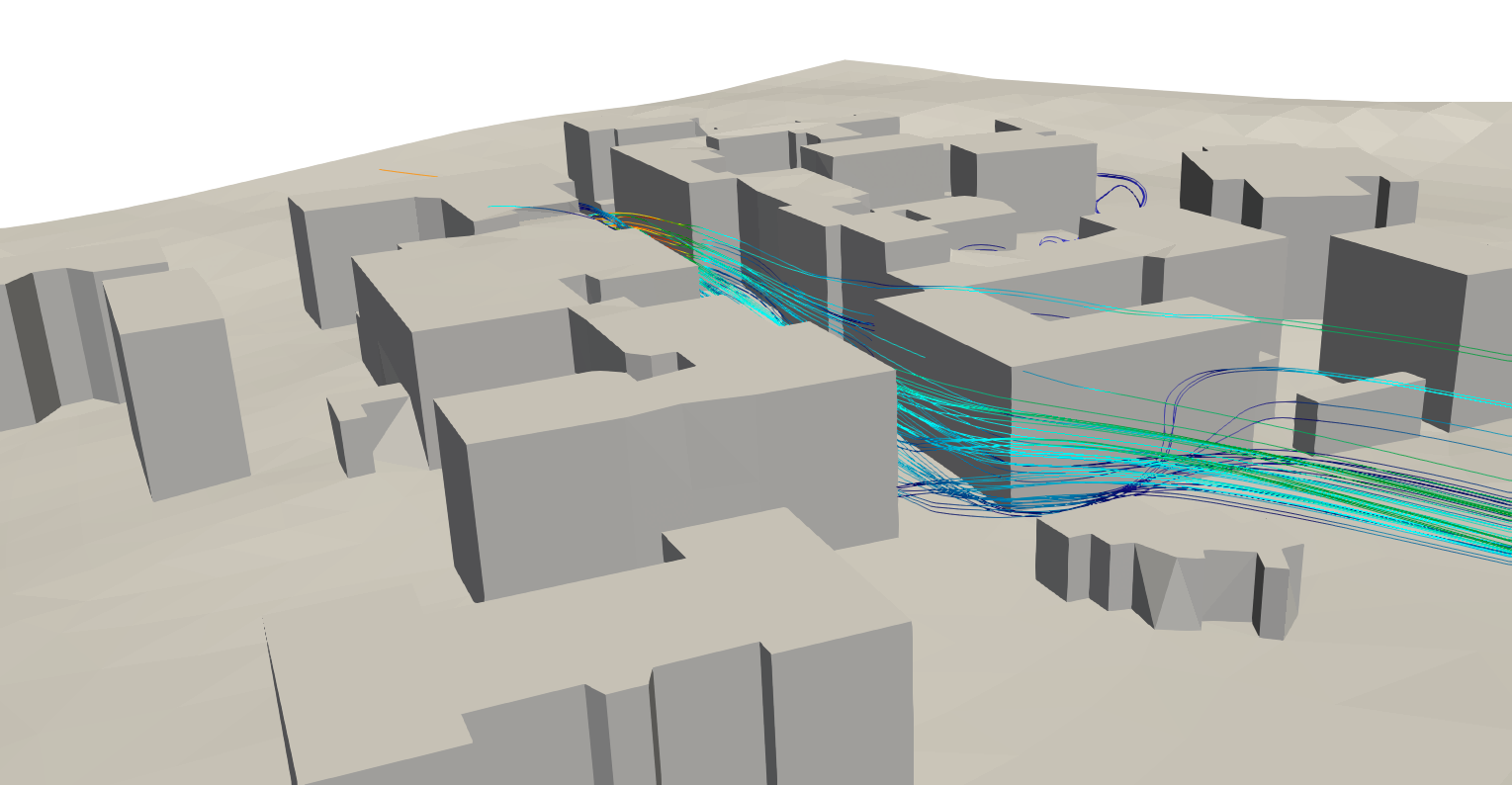}
    \caption{Streamlines of the velocity field.}
    \label{fig:streamlines}
\end{figure}

\begin{figure}[h]
    \begin{subfigure}[t]{0.45\textwidth}
        \centering
        \includegraphics[width=\textwidth]{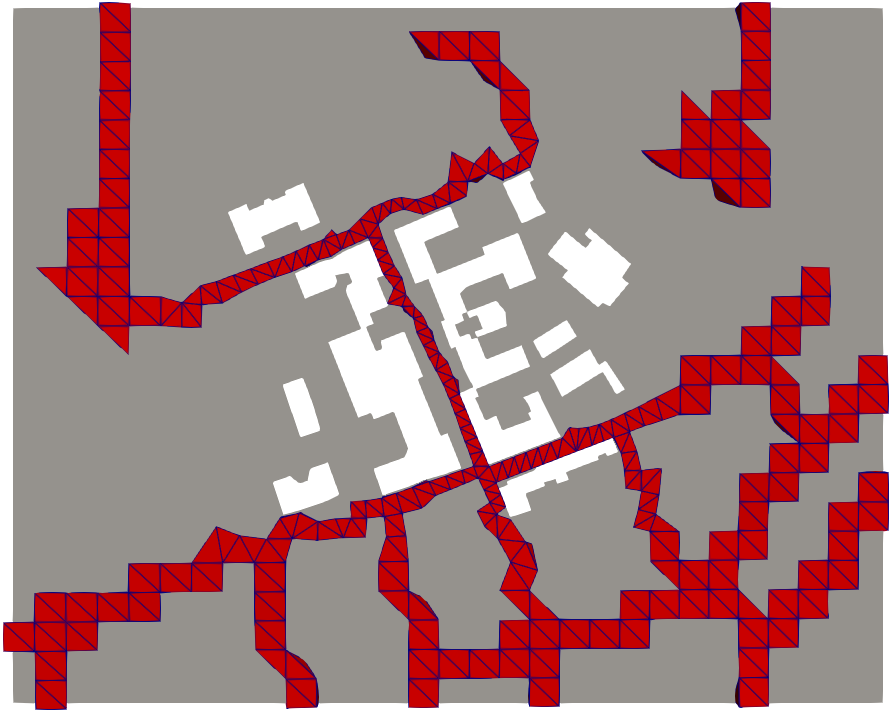}
        \caption{Red cells are those in which the source term is defined.}
        \label{fig:source_term_map-jpeg}
    \end{subfigure}
    \begin{subfigure}[t]{0.48\textwidth}
        \centering
        \includegraphics[trim=0 0 0 0,clip,width=\textwidth]{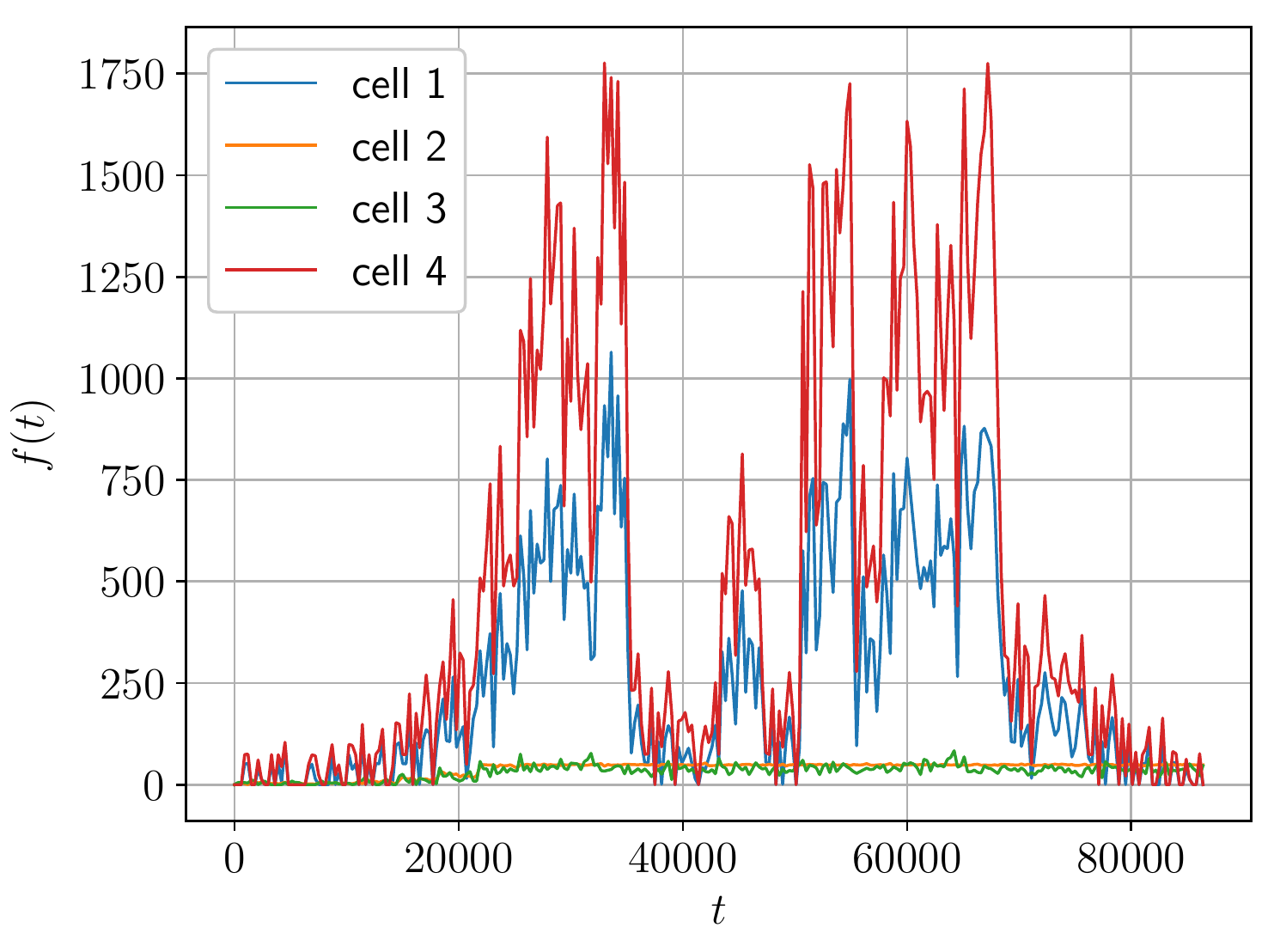}
        \caption{Pollution emission is treated as a per cell volumetric source that changes in time: example for 4 cells.}
        \label{fig:source_term_map-jpeg}
    \end{subfigure}
    \caption{Source term.}
    \label{fig:source-term}
\end{figure}

\label{sub:constant_convective_field}
The first test case is the situation where the convective field $\mathbf{u}$ is fixed and not parametric. Therefore, we are interested in the reduction of only the source term. This was done by exploiting two different strategies that still involve the extraction of a POD basis for both the concentration field and the emission source itself.
The convective field is obtained solving the steady-state $k$-${\epsilon}$ RANs equation by considering the following boundary conditions:
\begin{equation}\label{eq:navstokes-boundary}
    \begin{cases}
        \mathbf{u}  = \bm{0} &\mbox{ on } \Gamma_{FrontAndBack} \bigcup \Gamma_{Ground} , \\
        \mathbf{u} = \bm{(0,10,0)} &\mbox{ on } \Gamma_{In}  , \\
        (\nu\nabla \mathbf{u} - p\bm{I})\bm{n} = \bm{0} &\mbox{ on } \Gamma_{Out};
    \end{cases}
\end{equation}
where  $\nu=1.5e-5$, and the boundaries are shown in Figure~\ref{fig:boundaries}.
Figures~\ref{fig:velField} and \ref{fig:streamlines} show the velocity field and the streamlines of the velocity field, respectively computed by the solution of the steady-state RANs equations\footnote{The solution to the RANs system will be indicated with $\mathbf{u}$, instead of $\mathbf{\overline{u}}$, to ease the notation.} within the computational domain.
The source term $f$\footnote{In the following, we will omit the dependencies $(\boldsymbol{x},t)$ for the spatio-temporal fields $c$,$\mathbf{u}$, and $f$.}  has support $\Omega_{f} \subset \Omega$ and is defined through empirical series (see Figure \ref{fig:source-term}). In particular, we have considered time series in which the source is known at time instants $t_{k}= k h$, with $k\in \mathbb{N}$ and $h=300$. This situation models the case where such a source term is obtained through measurements with sampling frequency $\frac{1}{h}$. The source term in any other instant is obtained through a linear interpolation of the known values.
These time series are obtained using a simulated micro-traffic model through the open-source software SUMO \cite{Sumo}, whose results are subsequently converted into emission data using the COPERT emission model \cite{Copert}. Further information regarding the modeling of source emissions can be found in Appendix \ref{sec:source_emission_database}.

Let us then consider solving Equation \eqref{eq-transport}. In particular, we are interested in finding $c: \Omega \times[0, T) \to \mathbb{R} $, whereby taking $T=86400$, we are simulating an entire day.
The idea is to address the problem by applying the approach presented in Section \ref{sec:rom}. Specifically, at this stage, the only parameter is the time variable $t=\boldsymbol{\mu} \in [0, T)=\mathcal{P} $.
To this end, we used the snapshots coming from the FV resolution. Using the notation of Section \ref{sec:numerical_approximation}, we considered as $\mathcal{P}_\text{train}$ the discretized time interval $[0,20000]$, and $N_\text{train}=200$.

\begin{figure}[h]
    \begin{subfigure}[t]{0.45\textwidth}
        \centering
        \includegraphics[width=\textwidth]{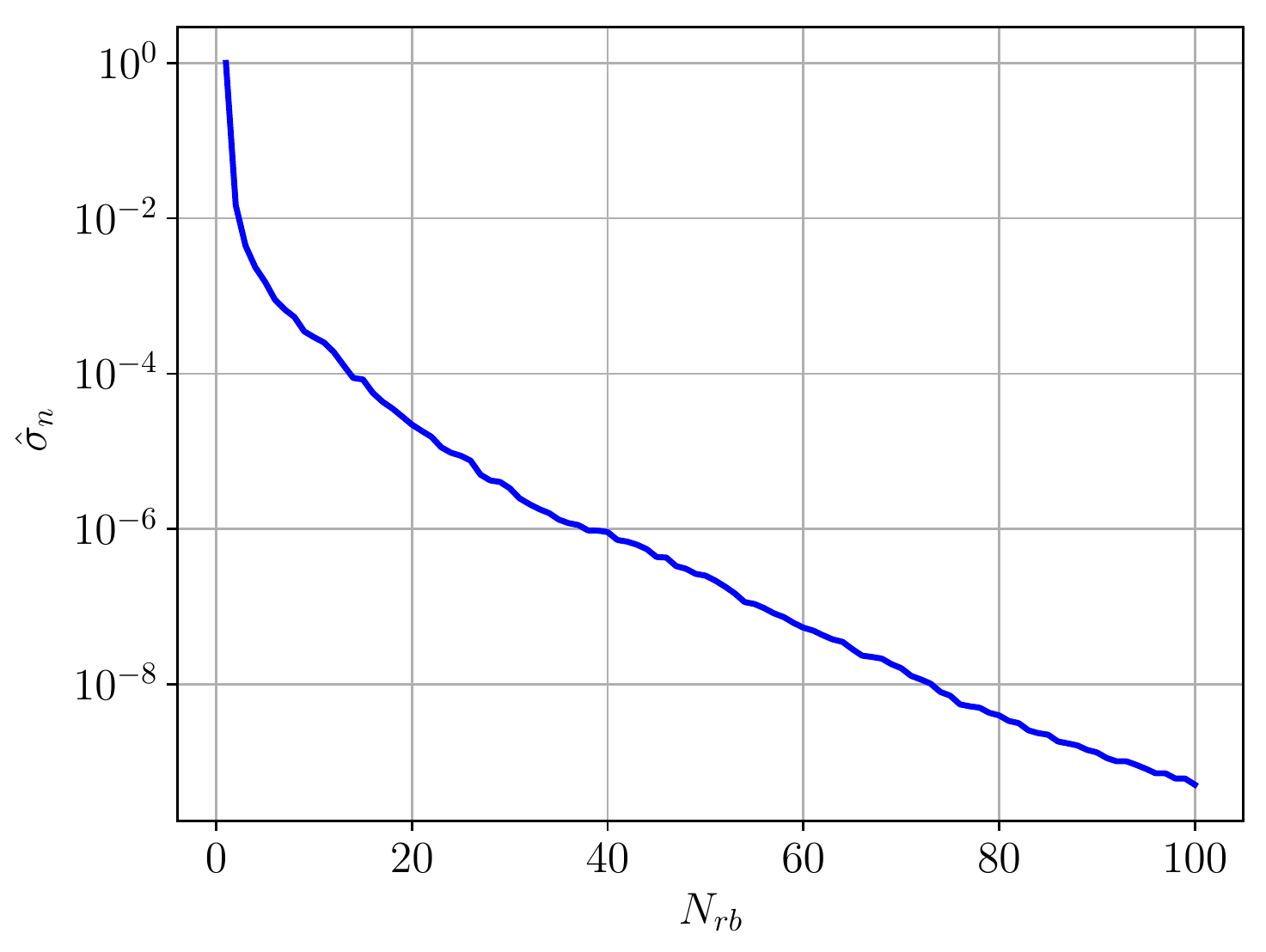}
        \caption{Decay of the normalized eigenvalues.}
    \end{subfigure}
    \begin{subfigure}[t]{0.45\textwidth}
        \centering
        \includegraphics[width=\textwidth]{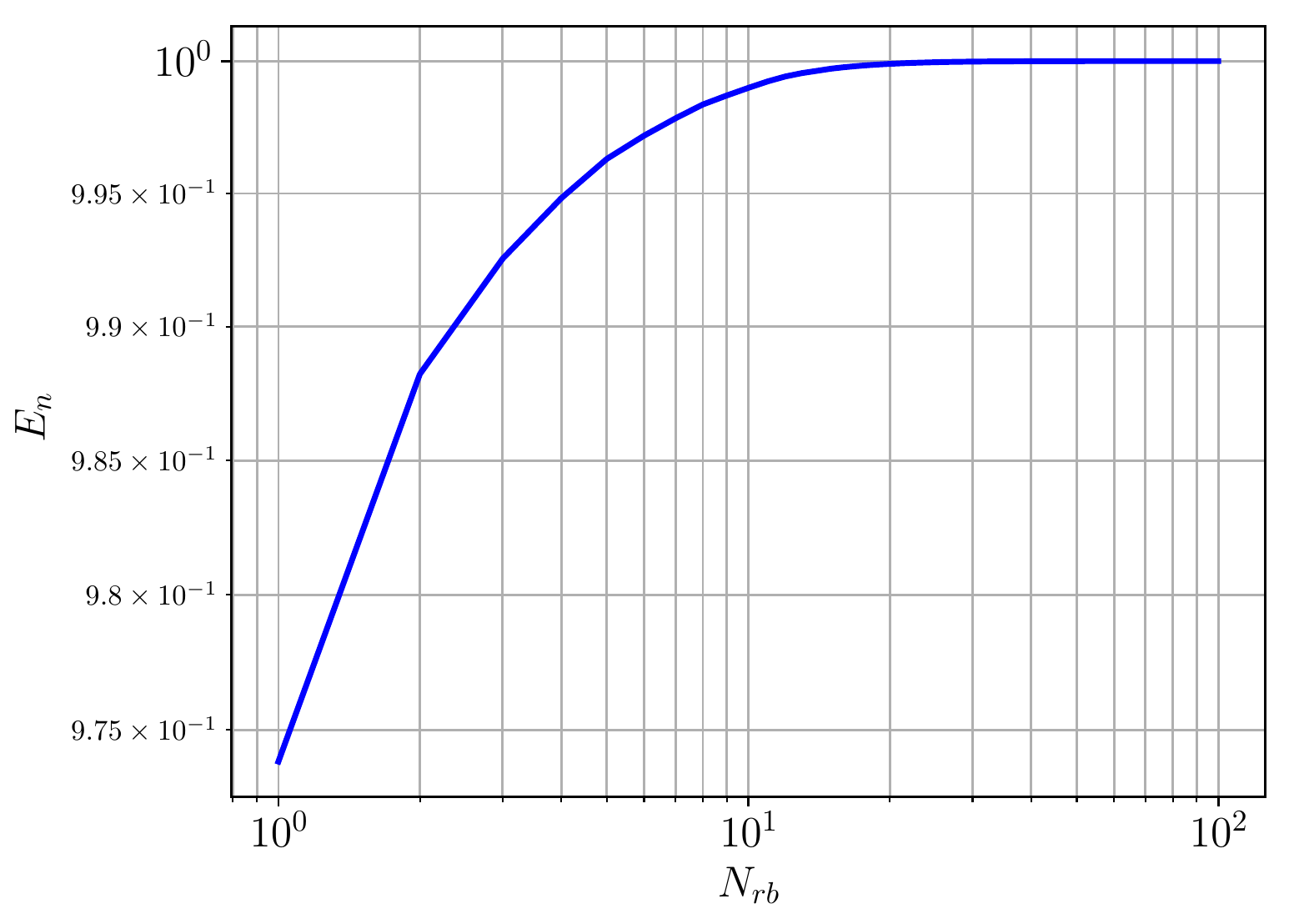}
        \caption{Retained energy.}
    \end{subfigure}
    \caption{Singular values indicators.}
    \label{fig:sing-values-indicators}
\end{figure}

\begin{figure}[h]
    \centering
    \includegraphics[width=0.5\textwidth]{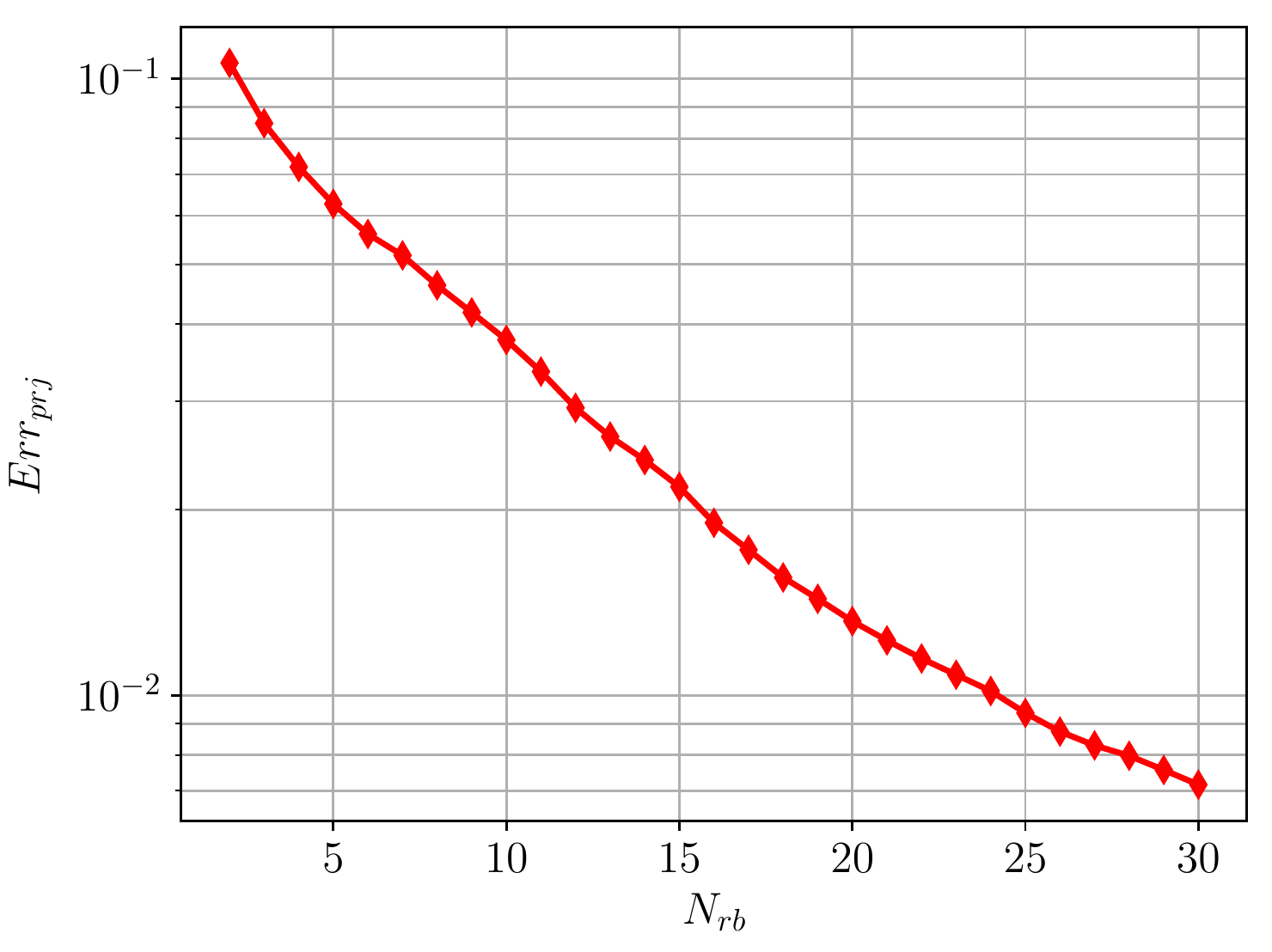}
    \caption{Average projection error in logarithmic scale.}
    \label{fig:avg-projection}
\end{figure}

The first indicator  to choose the cardinality of the reduced basis is the trend of the POD's singular values for the different fields. A general approach relies on a normalization w.r.t. the first eigenvalue, which is also dominant in magnitude:
\begin{equation}
    \quad \hat{\sigma_{i}}=\frac{\sigma_{i}}{\sigma_{1}}.
\end{equation}
Furthermore, if we take into account the cumulative sum of the eigenvalues, we get another important indicator known as \textit{retained energy}:
\begin{equation}
    E_{n}=\frac{\sum_{i=1}^{n} \sigma_{i}}{\sum_{i=1}^{N_\text{train}} \sigma_{i}}.
\end{equation}
We show in Figure \ref{fig:sing-values-indicators} the trend of both these indicators.
A valuable tool for evaluating the effectiveness of the basis extracted from the POD is the projection errors' trend as a function of the basis cardinality.
Since we are also interested in the use of the reduced model in a forecasting context, the projection error is taken as the average over the time interval $\mathcal{P}_\text{test}=[0,86400]$, consisting of $N_\text{test}=864$ instants of time, which includes the $200$ instants precedently used for the extraction of the POD basis:
\begin{equation}
    Err_{prj}=\frac{1}{N_\text{test}} \sum_{j=1}^{N_\text{test}} \frac{\| c_{h}(t_{j})-P_{n}(c_{h}(t_{j}))\|_{L^{2}(\Omega)}}{\|c_{h}(t_{j})\|_{L^{2}(\Omega)}};
\end{equation}
with $P_{n}$ being the projection operator of the FOM solution onto the reduced basis space.
The trend of this indicator is shown in Figure \ref{fig:avg-projection}. We observe an exponential decay of the error, which can be interpreted because the chosen POD basis is the one that minimizes the quantity in Equation \eqref{eq:minimize}. However, let us also recall that in this setting, we are using the reduced basis in the predictive regime to approximate temporal instances never seen by the model.

\begin{figure}[h]
    \begin{subfigure}[t]{0.45\textwidth}
        \centering
        \includegraphics[width=\textwidth]{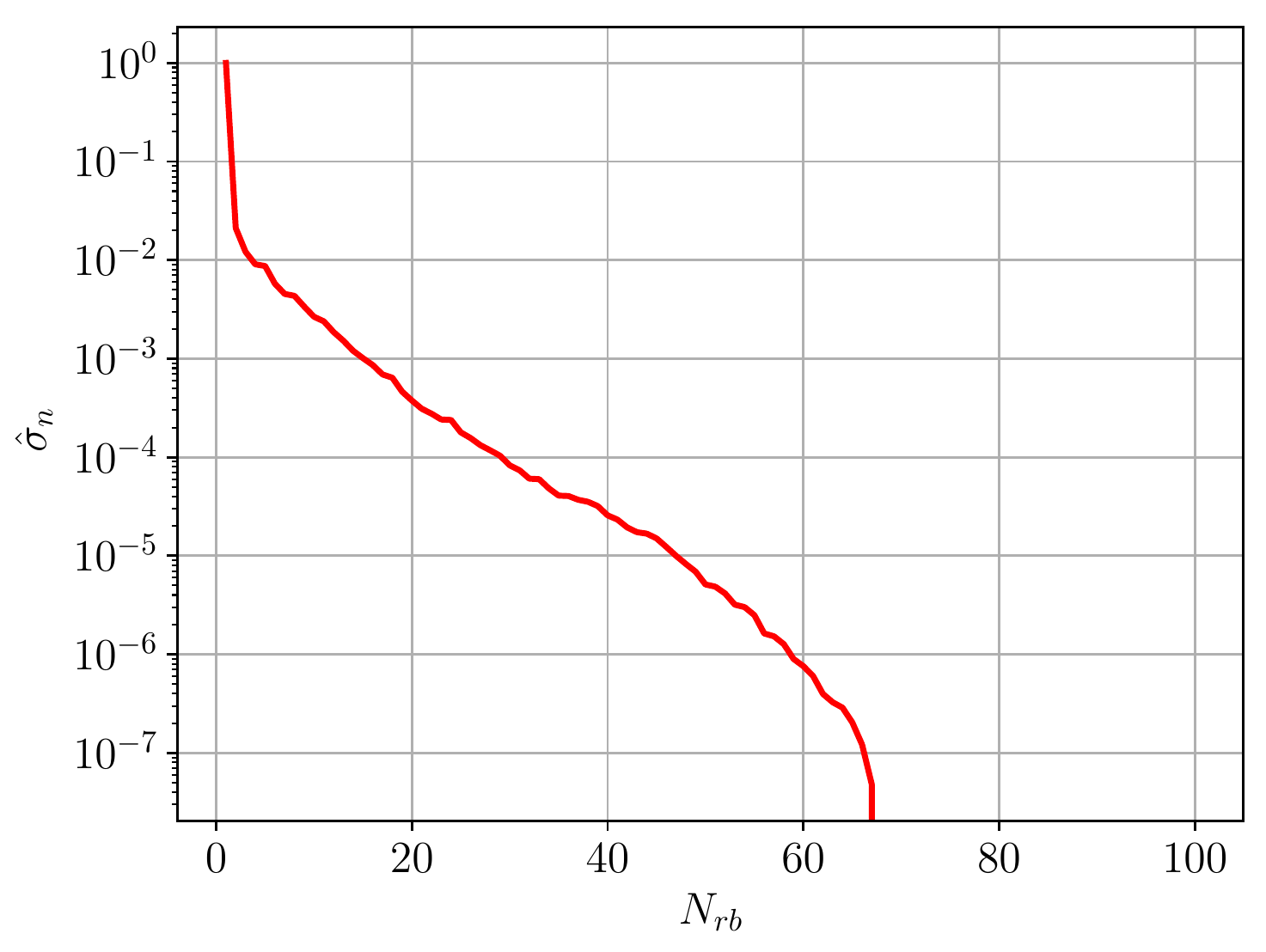}
        \caption{Decay of the normalized eigenvalues.}
    \end{subfigure}
    \begin{subfigure}[t]{0.45\textwidth}
        \centering
        \includegraphics[width=\textwidth]{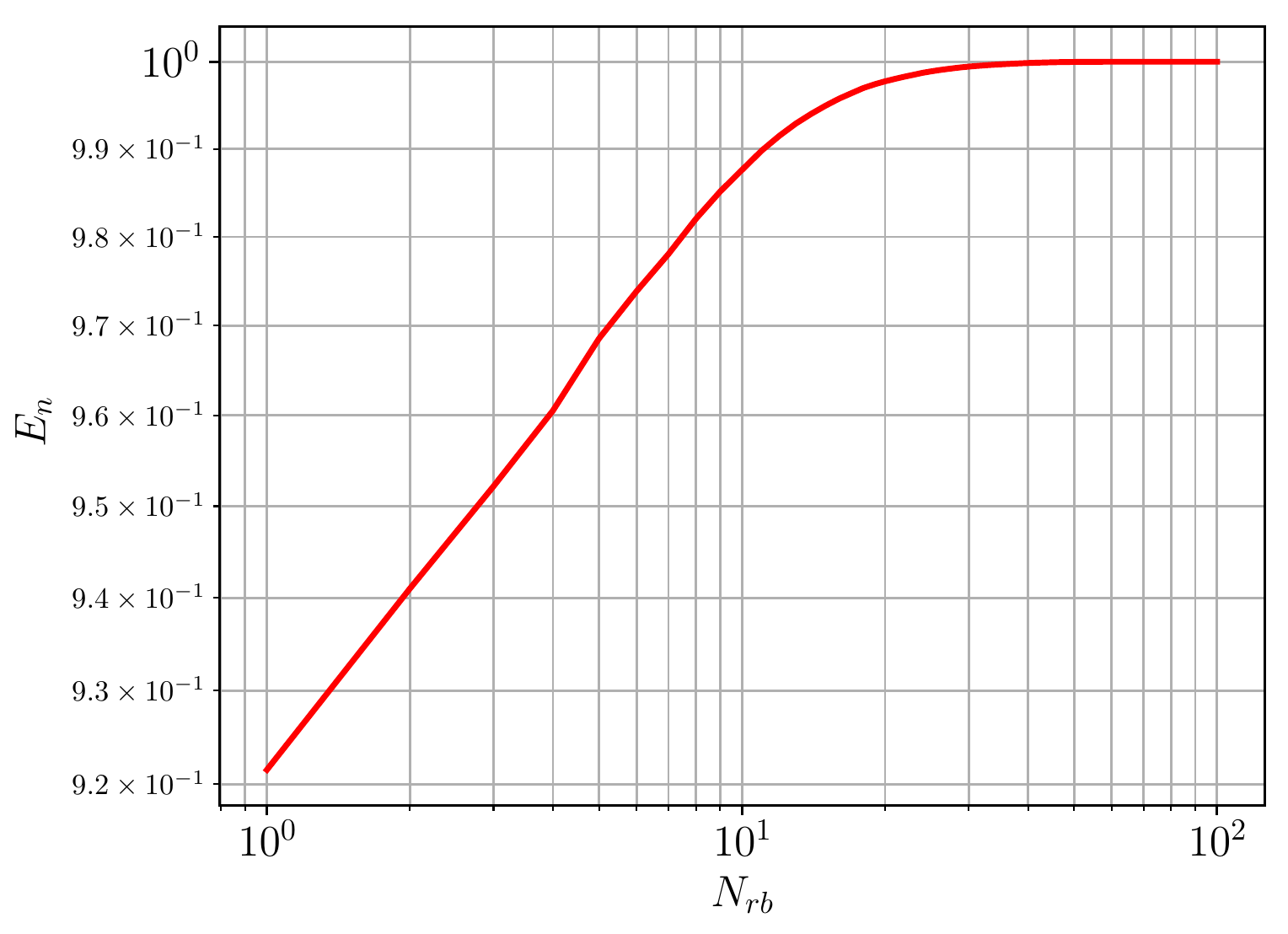}
        \caption{Retained energy.}
    \end{subfigure}
    \caption{Singular values indicators.}
    \label{fig:sing-values-indicators-S}
\end{figure}

\begin{figure}[h]
    \centering
    \includegraphics[width=0.8\textwidth]{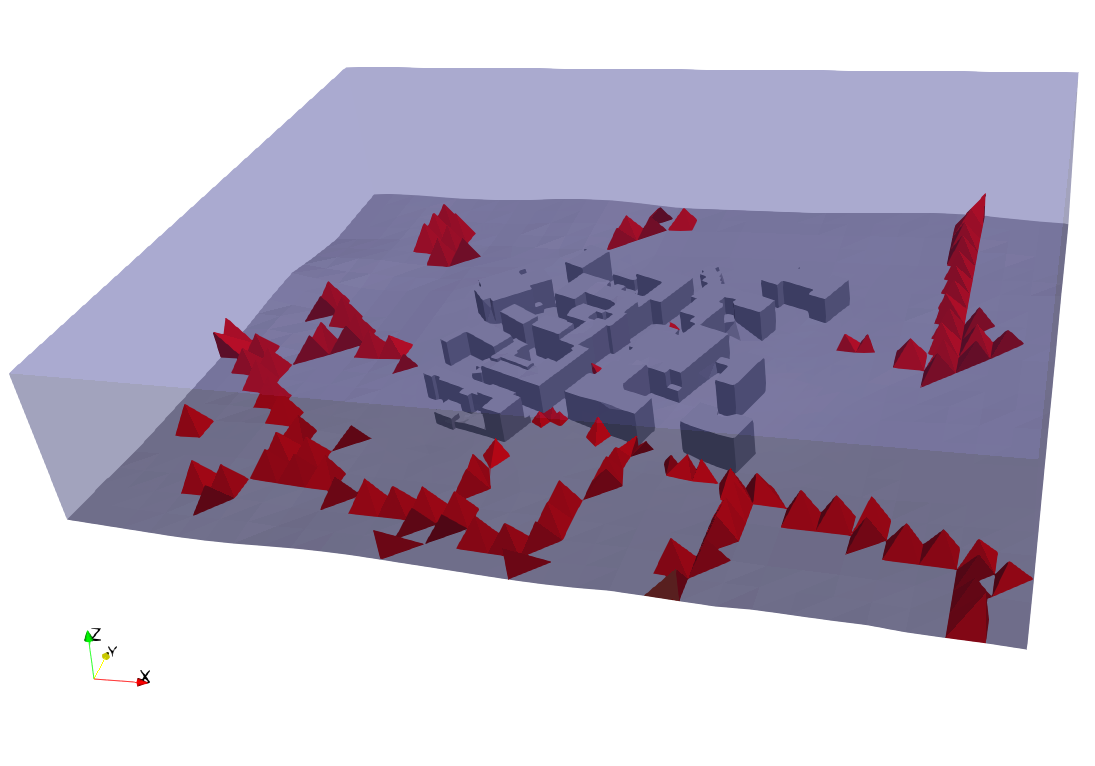}
    \caption{Magic points for $f(x,t)$  as computed from the DEIM algorithm.}
    \label{fig:magic-points}
\end{figure}

We investigated two possible approaches for the reduction of the source term.
The first strategy simply involves projecting the full order field $f_{h}$ onto the POD basis of $c$, i.e.:
\begin{equation}
    \label{eq:pod-proj-source}
    ( \boldsymbol{f}_{r}(t) )_{j} = (f_{h},\varphi_{j})_{L^{2}(\Omega)} \quad \text{for} \quad j=1, \ldots, N_{rb}.
\end{equation}
However, this strategy involves a computational cost that is still dependent on the FOM computational complexity because of the evaluation of the former integrals.
The second strategy investigated exploits hyper-reduction via the DEIM presented in Section \ref{sec:deim}.
In this case, extracting a POD basis for $f$ is necessary. In the same way, as for $c$, we report in Figure \ref{fig:sing-values-indicators-S} the POD indicators. The trend of these indicators is explained by the linear interpolation through which the field $f$ is obtained, resulting in a maximum number of linearly independent snapshots, which in this case are $66$.
Next are computed the indices for the magic points relative to the source term, presented in Figure \ref{fig:magic-points}. Since the magic points belong to a subset of the support of $f$, they are defined only on the roads of the computational domain.

\begin{figure}[h]
    \begin{subfigure}[t]{0.45\textwidth}
        \centering
        \includegraphics[width=\textwidth]{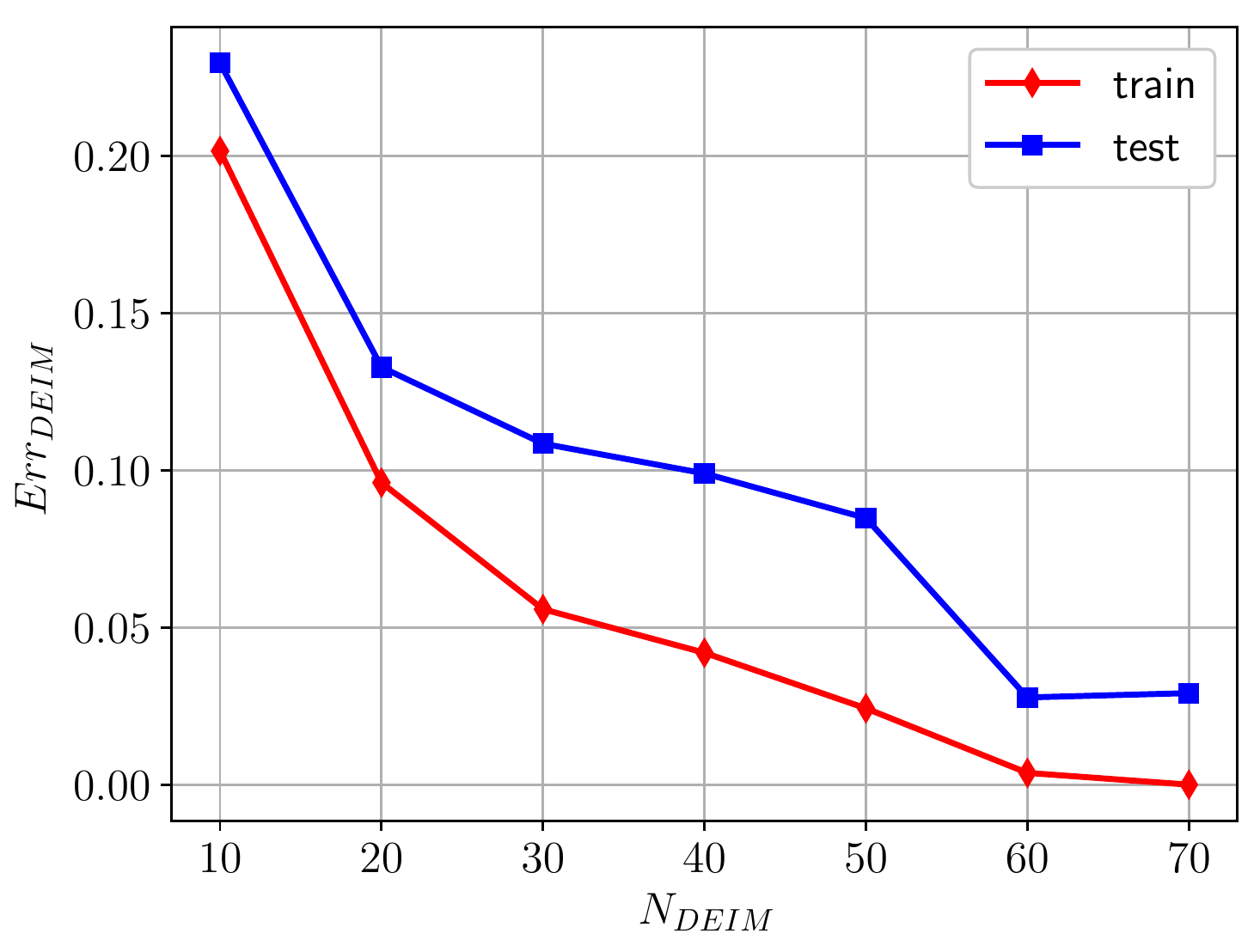}
        \caption{DEIM reconstruction error of the source $f$ on both the training and testing sets.}
        \label{fig:deim-error}
    \end{subfigure}
    \begin{subfigure}[t]{0.45\textwidth}
        \centering
        \includegraphics[width=\textwidth]{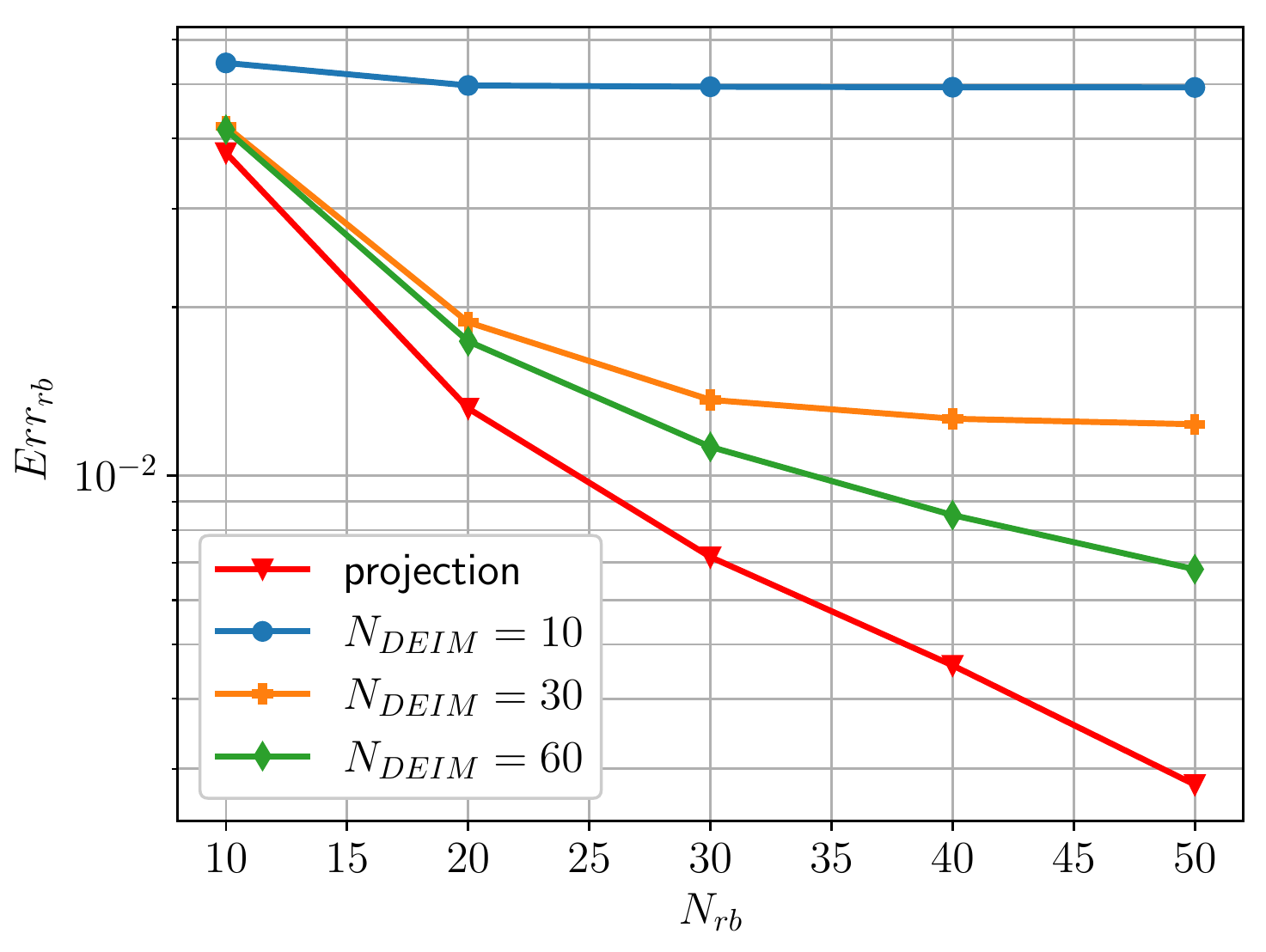}
        \caption{Average reconstruction error for $c$ in logarithmic scale.}
        \label{fig:rbm-error}
    \end{subfigure}
    \caption{}
\end{figure}

The average reconstruction error for the source term approximated using DEIM is presented in Figure \ref{fig:deim-error}.
In particular, we compare the mean reconstruction error evaluated on the set $\mathcal{P}_\text{train}=[0,20000]$ and the set $\mathcal{P}_\text{test}=[0,86400]$. We note an excellent reconstruction performance, which saturates when the number of DEIM bases exceeds the number of independent snapshots. In the latter case the error on the set $\mathcal{P}_\text{train}$ is zero while on $\mathcal{P}_\text{test}$ is equal to $\expnumber{2}{-2}$.
After obtaining the DEIM basis, it is possible to develop the ROM model.

We have analyzed the average error trend for some basis cardinalities in the solutions obtained through the RBM coupled with the DEIM reconstruction of the source term. In this case, the difference between the high-fidelity solutions and the reduced ones was measured as follows:
\begin{equation}
    Err_{rb}=\frac{1}{N_\text{test}} \sum_{j=1}^{N_\text{test}} \frac{\| c_{h}(t_{j})-c_{rb}(t_{j})\|_{L^{2}(\Omega)}}{\|c_{h}(t_{j})\|_{L^{2}(\Omega)}}.
\end{equation}

\begin{figure}[h]
    \centering
    \includegraphics[width=0.5\textwidth]{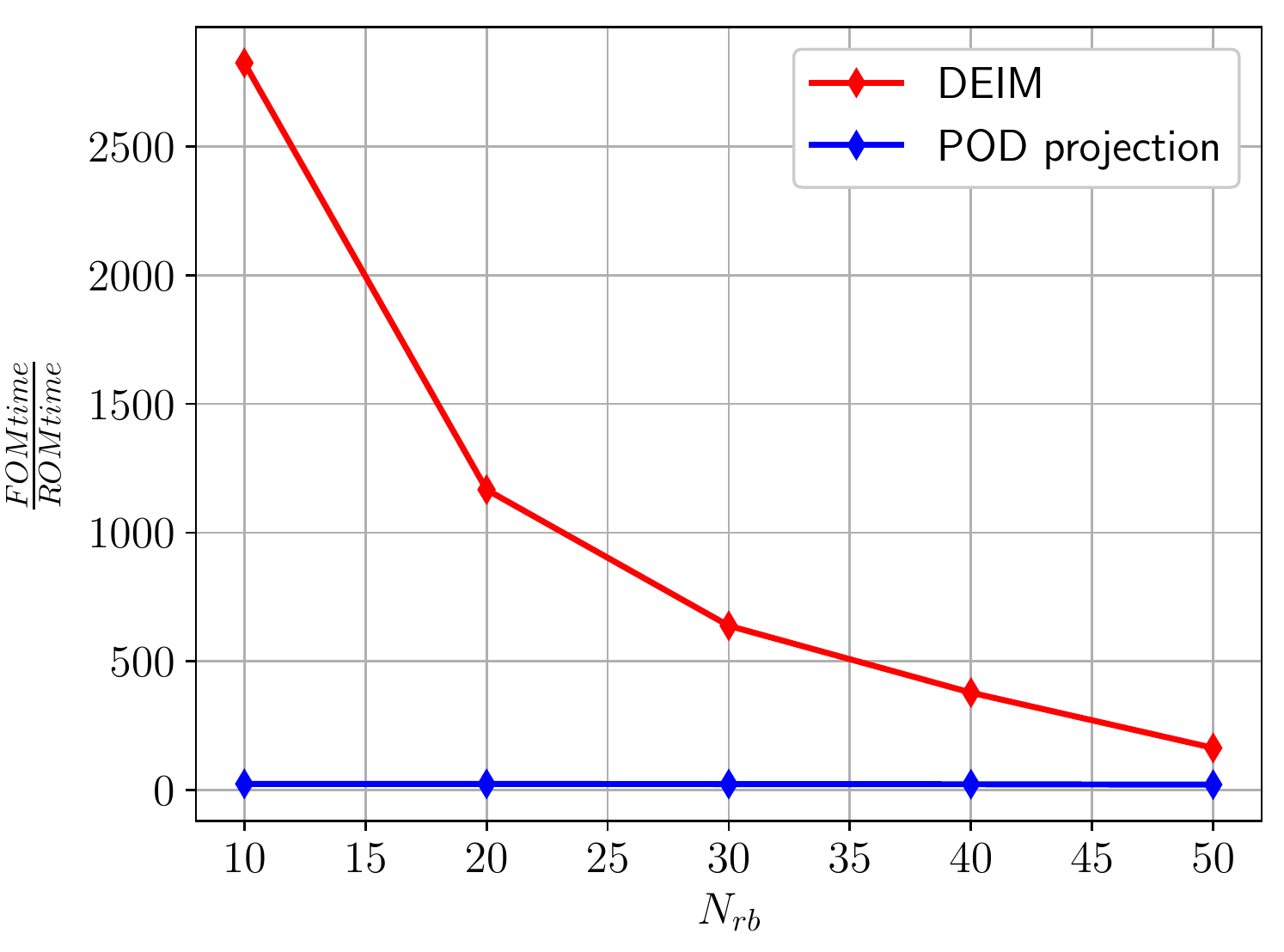}
    \caption{Speed-up comparison of DEIM model and POD projection.}
    \label{fig:plot_speedup-py}
\end{figure}

In Figure \ref{fig:rbm-error}, we observe that, as expected, the errors are higher than their direct projection counterparts.
In particular, we observe how increasing the number of bases for $c$ improves the model's accuracy. However, the lower bound on the error is always represented by how well we approximate the source through DEIM. Thus, we always arrive at a reduced basis cardinality, dependent on $N_\text{DEIM}$, which saturates the reconstruction error.
One might consider it appropriate to work with the maximum numbers $N_\text{DEIM}$ and $N_{rb}$. However, since for the present application, the real-time prediction is strictly related to the computational cost of solving the ROM model, one needs to analyze also the latter quantity. This analysis can be conducted by introducing the \textit{speed-up} of the ROM model, which is the quantity given by the ratio of the time required to solve the FOM problem to the time required to solve the ROM problem.
The trend of this quantity is shown in Figure \ref{fig:plot_speedup-py}, and one can observe its exponential decrease as a function of the number of reduced bases $N_{rb}$. The same figure also compares the speed-up obtained through DEIM with the alternative represented by the direct projection onto the POD basis (Equation \eqref{eq:pod-proj-source}). In the latter case, since the resolution implies the computation of integrals involving FOM fields, the speed-up is only about $24$ and does not scale with the number of bases used.

The previous analysis clarifies how the choice of cardinality $N_\text{DEIM}$ and $N_{rb}$ must be tailored to the particular application. The following results were obtained in the case where $N_\text{DEIM}=60$ and $N_{rb}=30$, obtaining a speed-up of approximately 640.

\begin{figure}[h]
    \centering
    \includegraphics[width=0.55\linewidth]{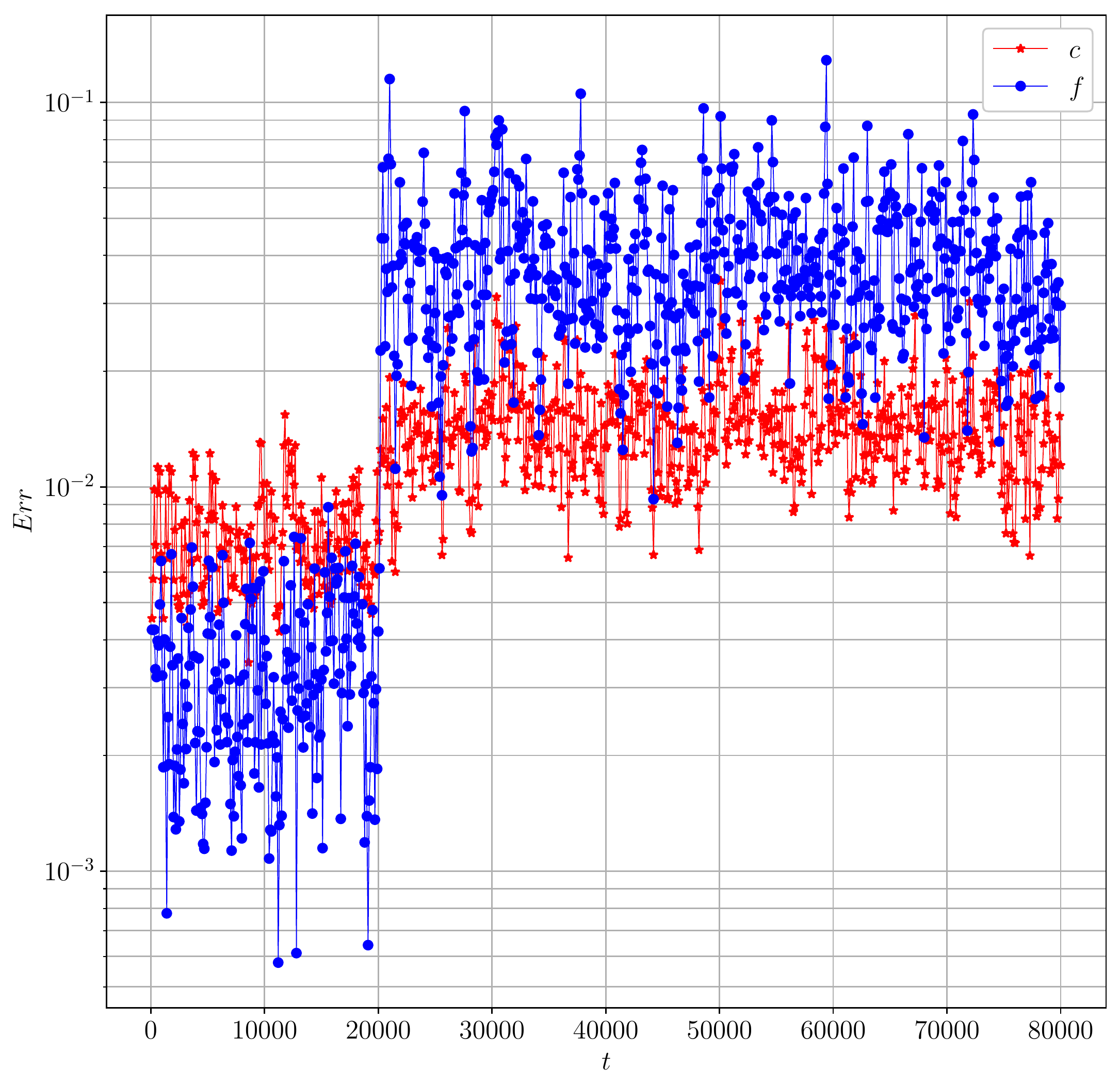}
    \caption{Relative reconstruction error in logarithmic scale with $N_{rb}=30$ and $N_\text{DEIM}=60$ over the time interval [0,80000].}%
    \label{fig:error_day}
\end{figure}

Following the implementation of the RBM, the reduced solutions for $\mathcal{P}_\text{test}$ are obtained. We then compared them with the full order counterpart through the relative errors for both $c$ and $f$ (Figure \ref{fig:error_day}).
We observe a highly oscillatory behavior of the errors as a function of time. In particular, the errors are lower on the time interval [0,20000] since the POD bases were extracted precisely from the snapshots belonging to it. However, surprisingly, despite a large increase in the source reconstruction error ($>10\%$) for $t>20000$, we witness a smaller increase for the $c$ error, which never exceeds $3.4\%$.

\begin{figure}[h!]
    \centering
    \includegraphics[width=1\linewidth]{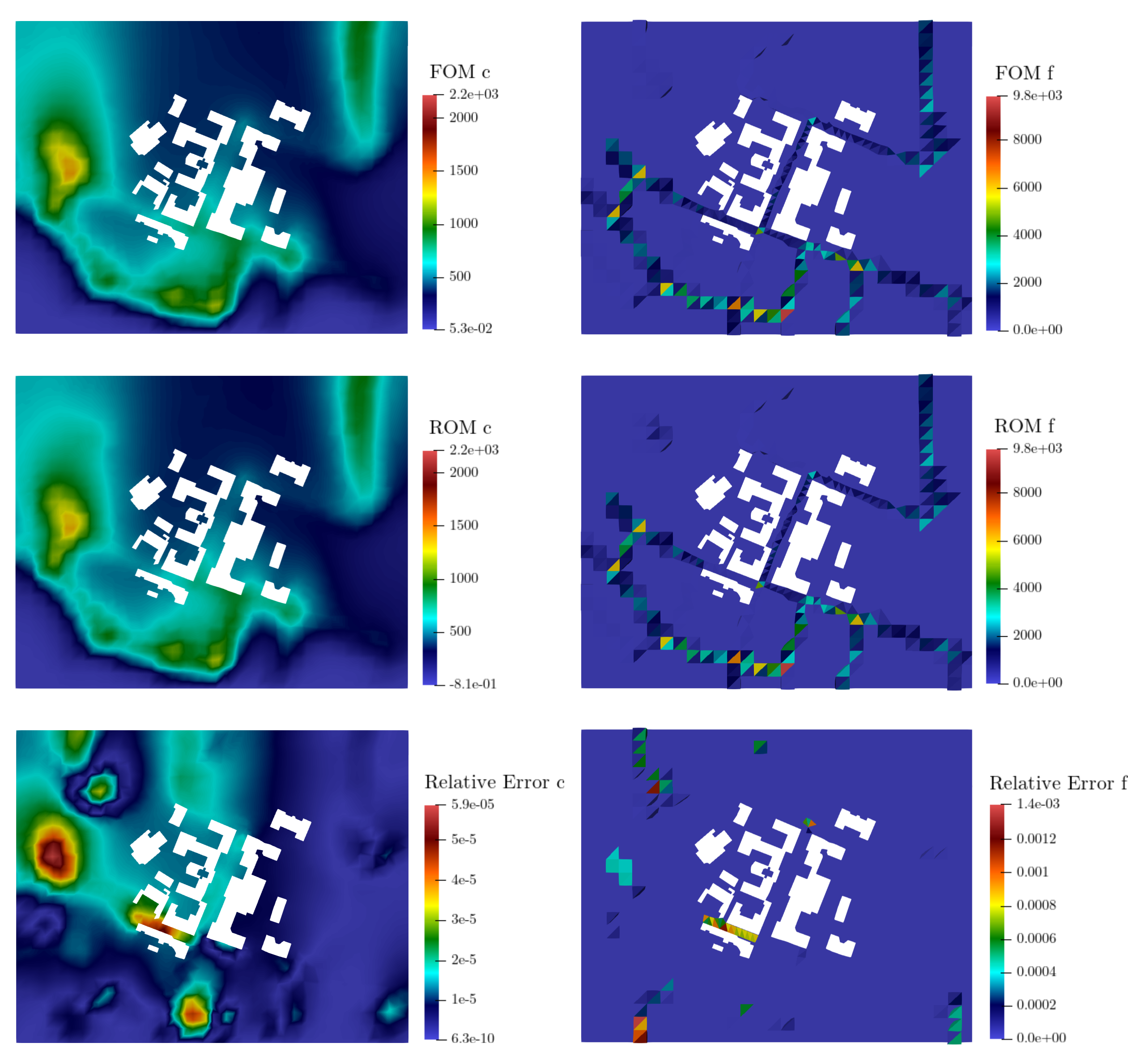}
    \caption{ Concentration and source emission fields for $t=50100$: FOM on top, ROM in the middle, and relative error between them on the bottom.}%
    \label{fig:Online-offline-deim-c}
\end{figure}

For some applications, it may also be important to understand the reconstruction's effectiveness as a function of the spatial coordinate $\boldsymbol{x}$. For this reason, we have also computed the relative error field, which assumes the following expression:
\begin{equation}
    \delta c_{rel} (\boldsymbol{x})=\frac{|c_{h}(\boldsymbol{x},t_{j})-c_{rb}(\boldsymbol{x,}t_{j})|}{{\|c_{h}(\boldsymbol{x},t_{j})\|_{L^{2}(\Omega)}}}, \quad    \delta f_{rel}(\boldsymbol{x})=\frac{|f_{h}(\boldsymbol{x},t_{j})-f_\text{DEIM}(\boldsymbol{x},t_{j})|}{{\|f_{h}(\boldsymbol{x},t_{j})\|_{L^{2}(\Omega)}}} .
\end{equation}
Figure \ref{fig:Online-offline-deim-c} shows the FOM, ROM, and error fields for both concentration and emission at time $t=50100$, which is the one with the higher global relative error ($Err_{rb}=0.0343$).
We can observe how the error for the concentration field traces the magnitude of the field itself, with a morphology similar to the error field for the emission.

\subsection{Parametric convective field}%
\label{sub:parametric_convective_field}

As we saw in the previous section, using RBM coupled with the DEIM hyper-reduction strategy allows us to effectively solve the linear advection-diffusion problem in the case where the emission field is defined through empirical time series.
However, to proceed toward a more accurate description of reality, the convective field $\mathbf{u}$ must be parameterized with respect to the boundary conditions.
In particular, we want to consider the following boundary conditions for the RANs equation:
\begin{equation}\label{eq:navstokes-param-boundary}
    \begin{cases}
        \mathbf{u}= \bm{0} &\mbox{ on } \Gamma_{Ground} , \\
        \mathbf{u} \cdot \bm{n}= \bm{0} &\mbox{ on } \Gamma_{Sky}  , \\
        \mathbf{u}= \mathbf{u}_{ref}\frac{\ln(\frac{z-d+z_0}{z_0})}{\ln (\frac{z_{r e f}+z_0}{z_0})} &\mbox{ on } \Gamma_{In} . \\
    \end{cases}
\end{equation}

\begin{figure}[hbt!]
    \centering
    \begin{subfigure}{0.45\textwidth}
        \centering
        \includegraphics[trim=80 50 100 100,clip,width=0.8\textwidth]{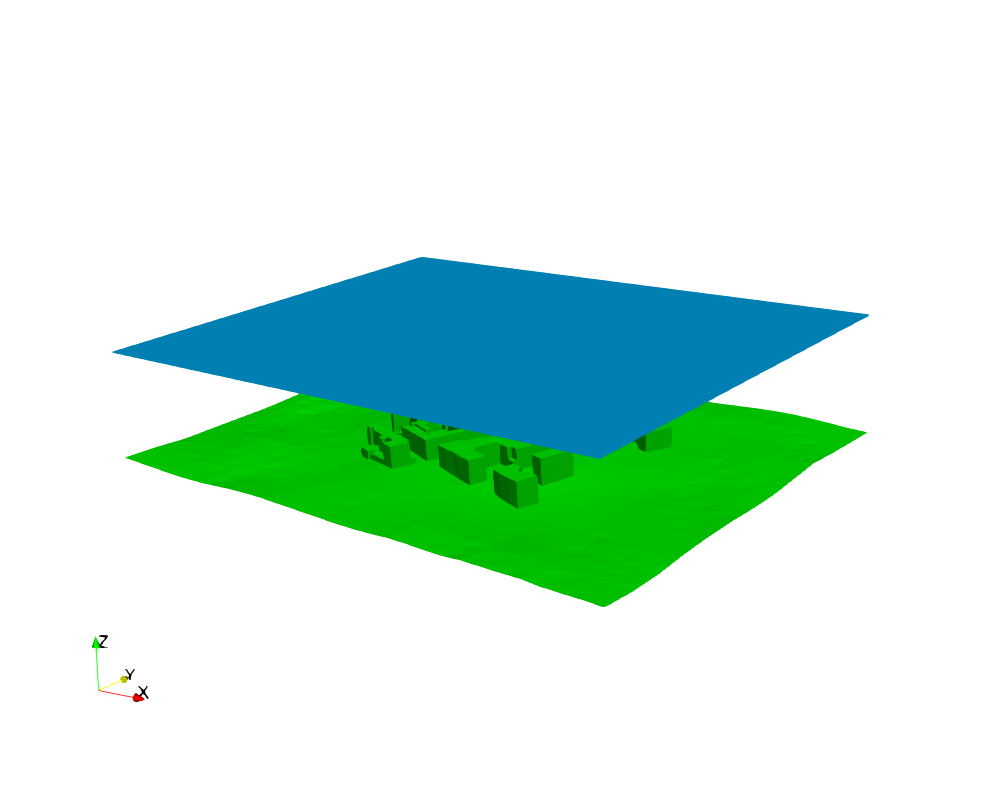}
        \caption{Sky boundary}
    \end{subfigure}
    \begin{subfigure}{0.45\textwidth}
        \centering
        \includegraphics[trim=80 50 100 100,clip,width=0.8\textwidth]{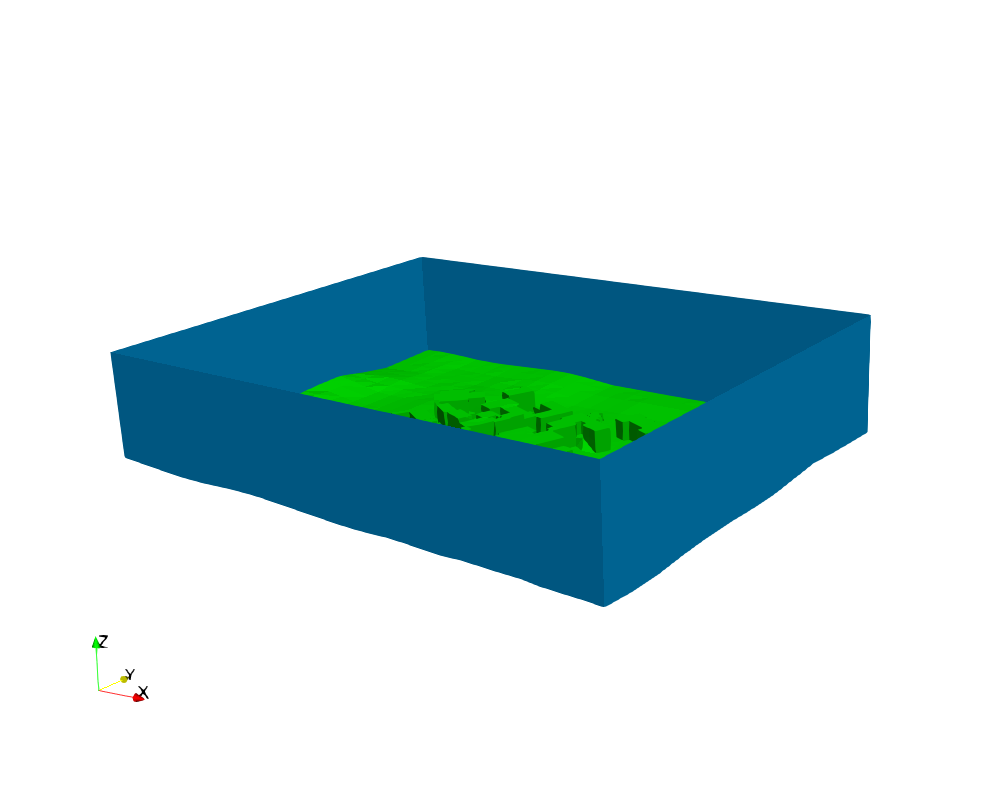}
        \caption{Inlet boundary}
    \end{subfigure}
    \caption{Boundaries of the computational domain, the ground boundary is reported in green.}
    \label{fig:boundaries-parametric}
\end{figure}
The boundary condition at the inlet $\Gamma_{In}$ (Figure \ref{fig:boundaries-parametric}) is  a log-law type ground-normal inflow boundary condition for wind velocity and the turbulence quantities for homogeneous, two-dimensional, dry-air, equilibrium and neutral atmospheric boundary layer (ABL) modelling \cite{Richards_1993}. The parameters for the ABL conditions are: the reference height $z_{ref}=10$, the ground-normal displacement height $d=0$ , and the aerodynamic roughness length $z_{0}=0.1$. In this setting, the reference velocity is parametrized with respect to its incidence and magnitude, such that $\mathbf{u}_{ref}= (\mu _{1}cos(\mu_2),\mu_{1} \sin (\mu_2),0)$ on $\Gamma_{In}$.

\begin{figure}[hbt!]
    \centering
    \begin{subfigure}[t]{0.49\textwidth}
        \centering
        \includegraphics[width=1\linewidth]{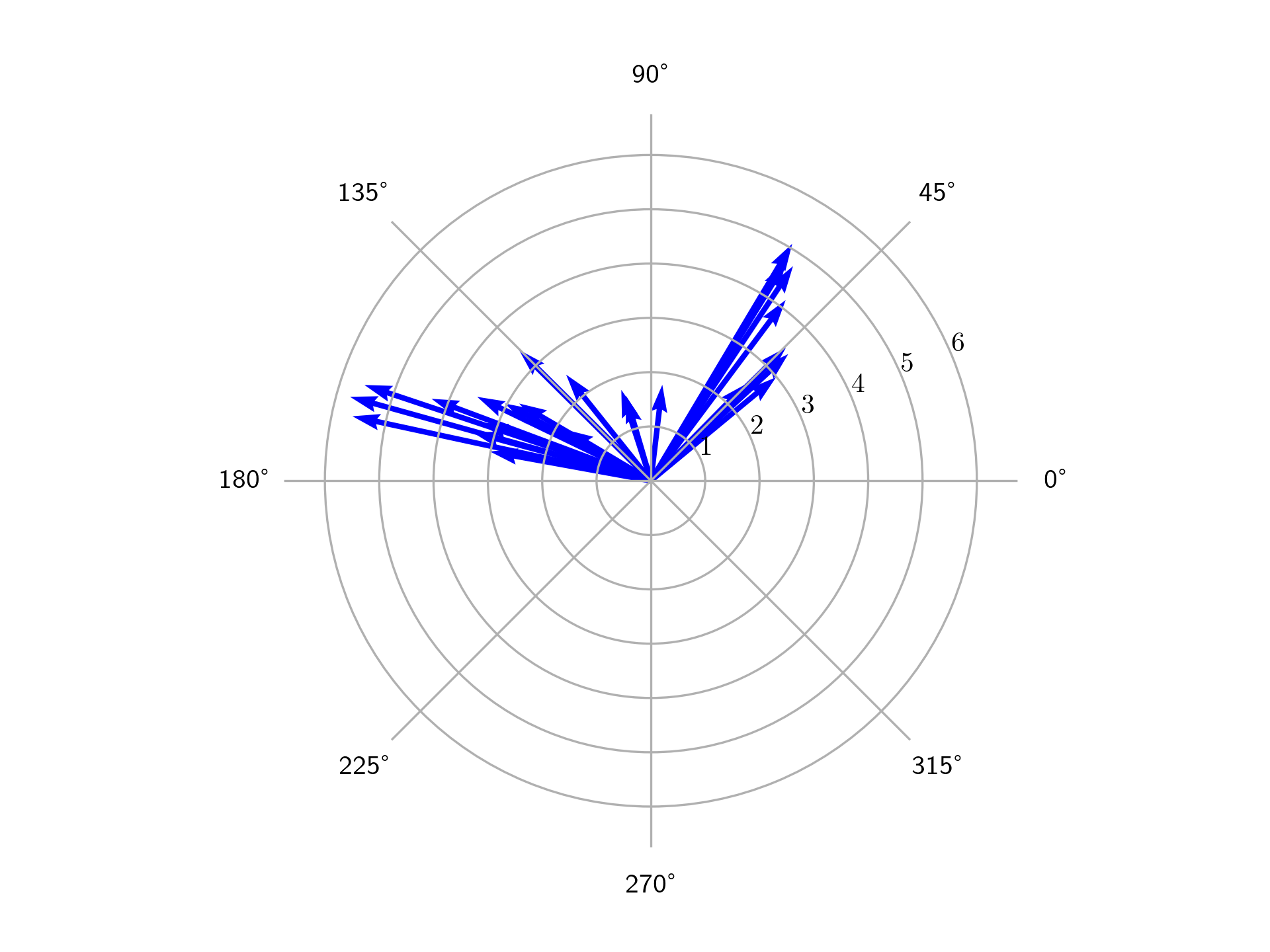}
        \caption{Wind rose for Day 0.}%
        \label{fig:wind_rose}
    \end{subfigure}
    \begin{subfigure}[t]{0.49\textwidth}
        \centering
        \includegraphics[width=1\linewidth]{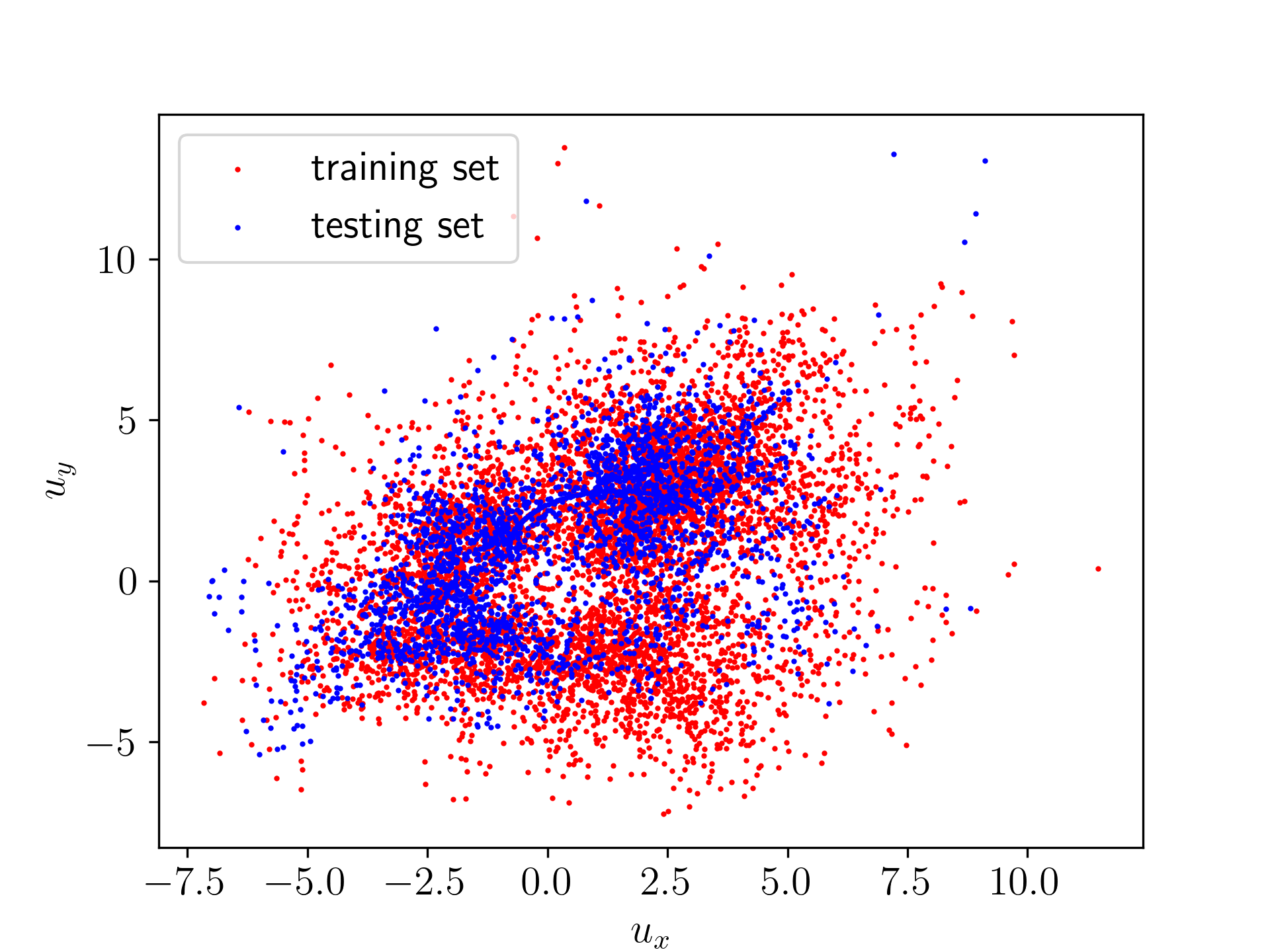}
        \caption{Scatter representation of inlet velocities for the simulated database.}%
        \label{fig:scatter_database}
    \end{subfigure}
    \caption{Inlet velocity conditions.}
\end{figure}

Therefore, we want to produce a ROM that can tackle the parametric convective field case, preserving the speed-up obtained in the previous test case. The motivation of the present investigation is to create a numerical framework for developing air quality monitoring models from a combination of information from a discrete number of local weather monitoring stations and urban emission data.
The first step for producing our model results in the application of the POD-R technique (Section \ref{sub:pod_with_inteporlation_podi_}) using a NN approach (Appendix \ref{sub:artificial_neural_networks}).
We decided to consider a database for weather conditions that included a large number of  revelations of wind conditions collected by a station data. Precisely, the database considered for inlet conditions consists of 7800 samples, corresponding to one detection every hour for 325 days in total.
We then divided the database in question, following a $70$-$30$  machine learning paradigm, into a dataset for training the model $\mathcal{D}^{\mu}_\text{train}$, and one for testing $\mathcal{D}^{\mu}_\text{test}$.
This subdivision corresponds to training the model with conditions related to $227$ training days and subsequently testing it on $98$ test days. An example of the inlet conditions for one of the simulated days is shown in Figure \ref{fig:wind_rose}, while Figure \ref{fig:scatter_database} presents the entire database.

We start by solving the full order steady RANs system for each of the inlet conditions belonging to $\mathcal{D}^{\mu}_\text{train}$. This allows us to obtain a POD basis for the velocity flow $\phi$.
The idea is to use the POD basis obtained through this step to train a multilayer feedforward neural network.
The neural network takes as input the two-dimensional parameter $\boldsymbol{\mu}=\left( \mu_{1}, \mu_{2}\right)$ corresponding to the reference velocity parametrization, and maps it to the POD coefficients of the flux field $\phi$ (Figure \ref{fig:nn-draw}). Specifically, the per-sample loss function used for training (Equation \eqref{nn-opt}) is the weighted squared Euclidean distance, with weights given by the singular values of the POD, to account for the different modal contributions:
\begin{equation}
    \mathcal{L}\left(\boldsymbol{\phi}_{true}, \boldsymbol{ \pi}_{NN}\right)=\sum_{j=1}^{N_{\phi}}(\lambda_{j}( \phi_{true,j}-\pi_{NN,j})^{2}) \ .
\end{equation}
\tikzset{>=latex} %
\colorlet{myred}{red!80!black}
\colorlet{myblue}{blue!80!black}
\colorlet{mygreen}{green!60!black}
\colorlet{mydarkred}{myred!40!black}
\colorlet{mydarkblue}{myblue!40!black}
\colorlet{mydarkgreen}{mygreen!40!black}
\tikzstyle{node}=[very thick,circle,draw=myblue,minimum size=22,inner sep=0.5,outer sep=0.6]
\tikzstyle{connect}=[->,thick,mydarkblue,shorten >=1]
\tikzset{ %
  node 1/.style={node,mydarkgreen,draw=mygreen,fill=mygreen!25},
  node 2/.style={node,mydarkblue,draw=myblue,fill=myblue!20},
  node 3/.style={node,mydarkred,draw=myred,fill=myred!20},
}
\def\nstyle{int(\lay<\Nnodlen?min(2,\lay):3)} %
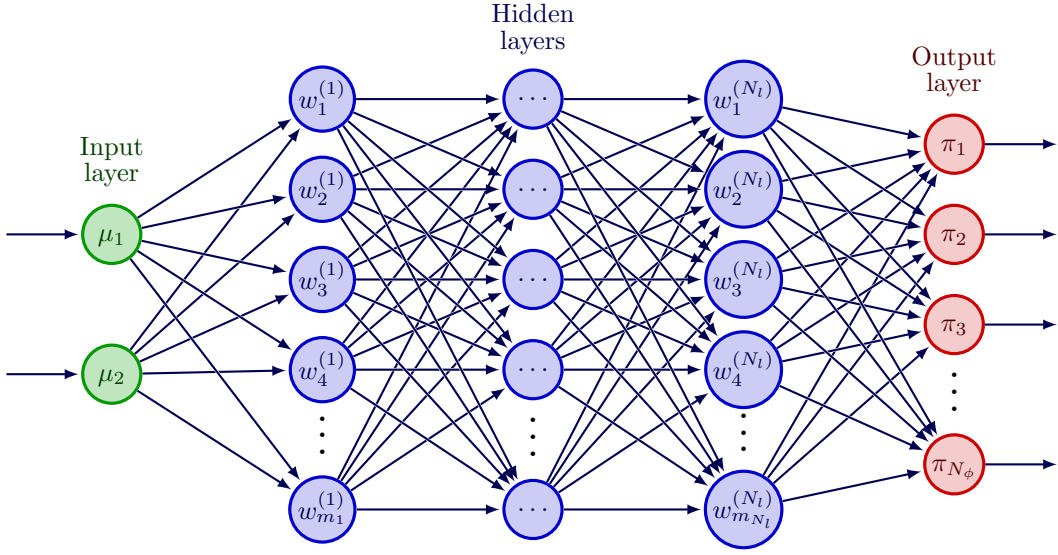
\begin{figure}
\centering
\begin{tikzpicture}[x=2.8cm,y=1.2cm]
  \readlist\Nnod{2,5,5,5,4} %
  \readlist\Nped{0,m_{1},0,m_{N_{l}},N_{\phi}} %
  \readlist\Nap{0,1,0,N_{l}} %
  \readlist\Nstr{2,m,N_{rb},4} %
  \readlist\Cstr{\mu,w^{(\Nap[\lay])},\pi} %
  \def\yshift{0.55} %

  \foreachitem \N \in \Nnod{
    \def\lay{\Ncnt} %
    \pgfmathsetmacro\prev{int(\Ncnt-1)} %
    \foreach \i [evaluate={\c=int(\i==\N); \y=\N/2-\i-\c*\yshift;
                 \x=\lay; \n=\nstyle;
                 \index=(\n>1 && \i==\N ?"\Nped[\lay]":int(\i));
                 }] in {1,...,\N}{ %
      \ifnumcomp{\lay}{=}{3} %
        {\node[node \n] (N\lay-\i) at (\x,\y) {$\ldots$} ;}
        {\node[node  \n] (N\lay-\i) at (\x,\y) {$\strut\Cstr[\n]_{\index}$};}

      \ifnumcomp{\lay}{>}{1}{ %
        \foreach \j in {1,...,\Nnod[\prev]}{ %
          \draw[white,line width=1.2,shorten >=1] (N\prev-\j) -- (N\lay-\i);
          \draw[connect] (N\prev-\j) -- (N\lay-\i);
        }
        \ifnum \lay=\Nnodlen
          \draw[connect] (N\lay-\i) --++ (0.5,0); %
        \fi
      }{
        \draw[connect] (0.5,\y) -- (N\lay-\i); %
      }

    }
    \ifnumcomp{\lay}{>}{1}
    {\path (N\lay-\N) --++ (0,1+\yshift) node[midway,scale=1.6] {$\vdots$};}{}
  }

  \node[above=3,align=center,mydarkgreen] at (N1-1.90) {Input\\[-0.2em]layer};
  \node[above=2,align=center,mydarkblue] at (N3-1.90) {Hidden\\[-0.2em]layers};
  \node[above=3,align=center,mydarkred] at (N\Nnodlen-1.90) {Output\\[-0.2em]layer};
\end{tikzpicture}
\caption{Scheme for a feedforward fully connected neural network that maps the parameter $\boldsymbol{\mu} \in \mathbb{R}^{2}$ to the flux POD coefficients $\boldsymbol{\pi_{NN}} \in \mathbb{R}^{N_{\phi}}$.}
\label{fig:nn-draw}
\end{figure}

\begin{figure}[h]
    \centering
    \includegraphics[width=0.6\linewidth]{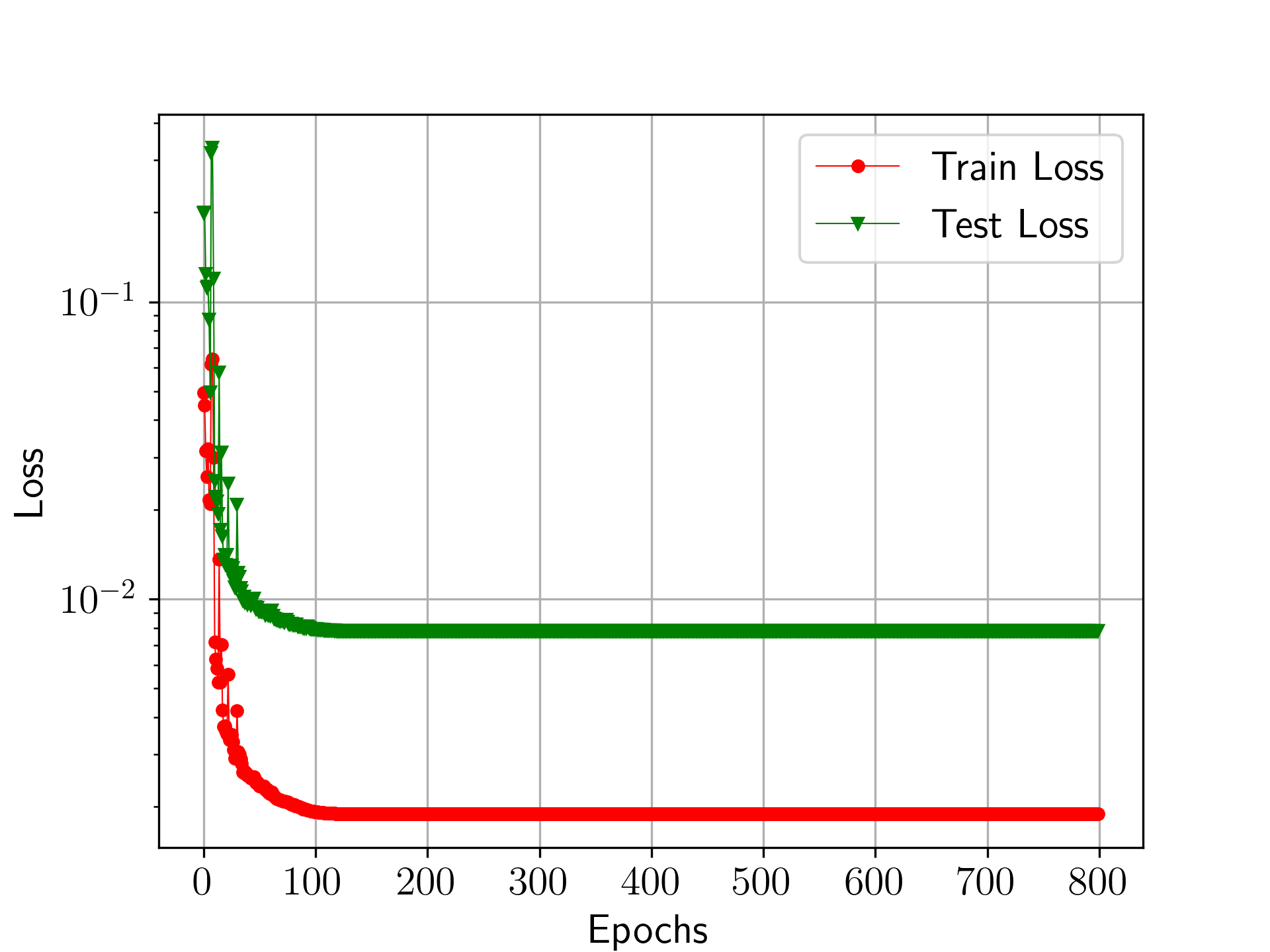}
    \caption{Loss function decay for both training and test sets.}%
    \label{fig:Loss}
\end{figure}
It should be noted that the generic target output is obtained in this case via a direct projection of the training full order fluxes on the POD basis.  The characteristics of the neural network are summarized in Table \ref{tab:table-nn}. For a more detailed discussion of the meaning of the various terms we refer the interested reader to  \cite{goodfellow}.
Following the optimization procedure (Figure \ref{fig:Loss}), the neural network can be used in the predictive regime.
\begin{table}[h]
    \centering
    \begin{tabular}{ |l|c| }
        \hline
        Typology &\textit{ Feedforward fully connected }\\ \hline
        Number of hidden Layers &\textit{ 6 }\\ \hline
        Input Neurons &\textit{ 2 }\\ \hline
        Output Neurons &\textit{ 200 }\\ \hline
        Activation function &\textit{ ReLu }\\ \hline
        Optimizer &\textit{ Adam with minibatching}\\ \hline
        Learning rate &\textit{1e-3}\\ \hline
        Weight decay &\textit{1e-8}\\ \hline
        Epochs  &\textit{80000}\\
        \hline
    \end{tabular}
    \caption{Neural network characteristics.}
    \label{tab:table-nn}
\end{table}

The second step involves the extraction of a POD basis for the concentration field $c$ and the source term $f$.
The procedure is similar to that presented in Section \ref{sub:constant_convective_field}, the only difference being it involves multiple simulation days. Specifically, we have 93600 data for emissions, corresponding to a detection every 300 seconds for 325 total days (see Appendix \ref{sec:source_emission_database}).
However,  we also simplified the source modelization by considering its support limited to cells belonging to the main roads of the computational domain, for a total of 169 cells.
Similar to what was done for the velocity conditions at the inlet, these data are divided in turn into a database for training $\mathcal{D}^{f}_\text{train}$, and one for testing $\mathcal{D}^{f}_\text{test}$.
Specifically, in this case, the training dataset contains the emissions for the time interval [0,10200] of each of the days represented in $\mathcal{D}^{\mu}_\text{train}$.

We note that now $\mathcal{P}\subset \mathbb{R}^{3}$, because in addition to the two parameters for the inlet condition, we also consider time $t$. Following the FV resolution we extract a POD basis for $c$ and $f$, by considering snapshots belonging to $\mathcal{P}_\text{train}=\{(\mu_{1}(t), \mu_{2}\left( t \right) ,t) | t \in \mathcal{T}_\text{train}\}$, where the set $\mathcal{T}_\text{train}$
consists of the union of the discrete intervals containing the first 10200 seconds of each day of the training set, for a total of $N_\text{train}$=7718 samples. We also point out that the functions $(\mu_{1}(t),\mu_{2}(t))$ are obtained through a linear interpolation of the discrete detections contained in $\mathcal{D}^{\mu}_\text{train}$.

Following the extraction of the above bases, we proceed with the production of the ROM model, according to the procedure given in Sections \ref{sec:rom}-\ref{sub:podnn-g}.
This means that the output of the neural network provides $online$ the vector of coefficients $\boldsymbol{ \pi}_{NN}$ that is used to obtain the matrix for the convective term from the tensor $\mathbf{\Gamma}$ as per Equation \eqref{eq:tensor-contraction}. In contrast, the source is still reduced using DEIM.

\begin{figure}[h]
            \centering
            \includegraphics[width=0.7\linewidth]{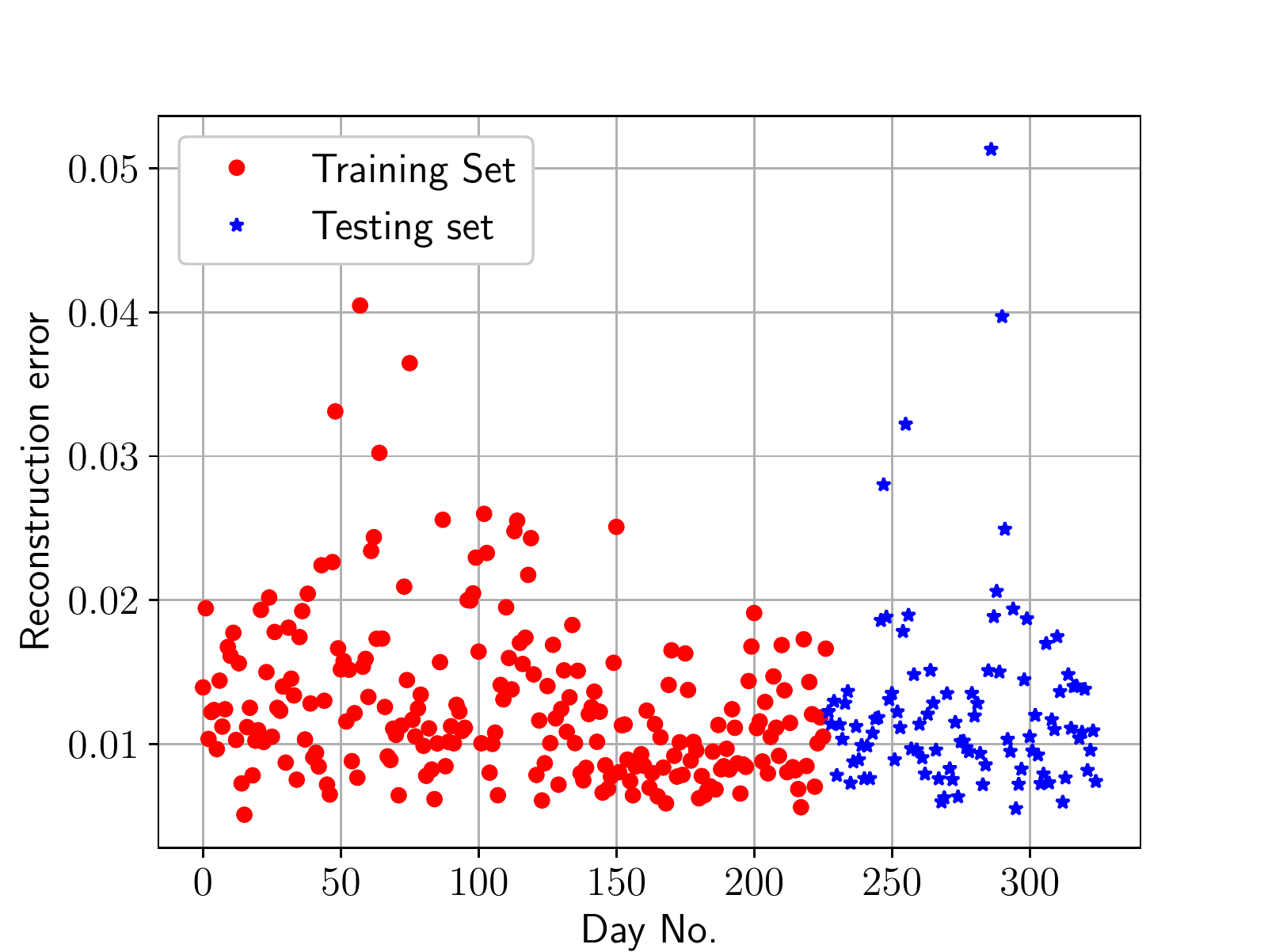}
            \caption{Average daily reconstruction error for Model D (Table \ref{tab:table-podgnn}).}%
            \label{fig:avg_error}
        \end{figure}

It should be noted that at this point, consistently with the considerations already stated in Section \ref{sub:constant_convective_field}, a choice must be made regarding the cardinality of the reduced bases of $\phi$,
$f$ and $c$, which we recall are $N_{\phi}, N_\text{DEIM}$ and $N_{rb}$. This choice affects both the computational cost of producing the ROM and the final accuracy of the predictions.
Regarding the first of the two aspects, we highlight how precisely the calculation of the convective tensor $\mathbf{\Gamma}$ is the bottleneck of the process.
We explored different size choices for the reduced model, which we tested on novel instances of both the convective field and the emission term.
In particular, we simulated the 325 days for which we had velocity boundary conditions and empirical time series for the emission.
Table \ref{tab:table-podgnn} shows the average relative error and speed-up for different ROMs obtained by changing the number of bases.
In particular, we can see that a low $N_{rb}$ number (Model A) leads to low accuracy but high speed-up. Increasing $N_{rb}$, when both $N_\text{DEIM}$ and $N_{\phi}$ are low, leads to decreades computional speed-up without an improvement in the model accuracy (Model B).
It is evident that increasing the cardinality of both the discrete empirical interpolation method (DEIM) basis functions, $N_\text{DEIM}$, and the number of reduced basis $N_{\phi}$, improves the accuracy of the reconstructed convective field and source term. This, in turn, can enhance the overall accuracy of the model (Models C and E). However, this behavior becomes saturated for high cardinalities, as observed in Model D. Moreover, the computational speed-up depends mainly on the cardinality of the reduced basis functions $N_{rb}$ (Models B, C, and D).
The results for the average daily relative error of Model D for the concentration field $c$ are shown in Figure \ref{fig:avg_error}.

        \begin{table}[h]
            \centering
            \begin{tabular}{ |c|c|c|c|c|c| }
                \hline
                \textbf{Model} & $N_{\phi}$ & $N_\text{DEIM}$ & $N_{rb}$ & Mean relative error & Speed-up  \\ \hline
                A & 10 & 10 & 10 & 0.101258 & 86404 \\
                B & 10 & 10 & 50 & 0.101028 & 5630 \\
                C & 20 & 20 & 50 & 0.012637 & 5913 \\
                D & 50 & 50 & 50 & 0.012630 & 5516 \\
                E & 20 & 20 & 70 & 0.010718 & 3009 \\
                \hline
            \end{tabular}
            \caption{Model statistics for different choices of the basis cardinality.}
            \label{tab:table-podgnn}
        \end{table}

\section{Conclusions and perspectives}%
\label{sec:conclusions_and_perspectives}

This paper presents a novel ROM model for the linear-advection diffusion problem in an urban setting. We proposed a new methodology that exploits a mixture of intrusive and nonintrusive model reduction, achieving to maintain the advantages of both approaches. The technical difficulty of the problem addressed lies in the nature of the convective field and the source term. The latter term was modeled through empirical time series to emulate the model's use in conjunction with real emission data. In the first test case, we considered a nonparametric convective field. In that case, we demonstrated the flexibility of the discrete empirical interpolation strategy, which allowed a complete decoupling between offline and online phases and a consequent significant computational speed-up. We employed a feedforward neural network to address the test case with a parametric convective field. This allowed us to reduce the phase related to solving the RANs equations, and integrate it directly with the Galerkin projection of the transport equation, thus preserving the computational advantages of the nonparametric case.

During the conclusion of this work, we have identified many interesting topics to be deepened.
First, our work was limited to the analysis of the linear advection-diffusion equation. A natural extension corresponds to using more accurate air pollution models, which also consider chemical transformation and ground deposition. Typically, the more sophisticated air quality models present additional nonlinearities. However, these can be addressed like what has already been proposed in this paper. In addition, the nonintrusive nature of our approach also allows us to consider a different CFD urban model. Thus, we could explore the use of a LES approach for the convective term.

Finally, an interesting perspective is the integration of our framework with an inverse modeling approach, which would allow the geometric parameters of the computational domain to be optimized as a function of local or global indicators of pollutant concentration.

\section{Acknowledgments}
This work was partially funded by European Union Funding for Research and Innovation — Horizon 2020 Program — in the framework of European Research Council Executive Agency: H2020 ERC CoG 2015 AROMA-CFD project 681447 "Advanced Reduced Order Methods with Applications in Computational Fluid Dynamics" P.I. Professor Gianluigi Rozza. We also acknowledge the PRIN 2017 "Numerical Analysis for Full and Reduced Order Methods for the efficient and accurate solution of complex systems governed by Partial Differential Equations" (NA-FROM-PDEs).
The research for this paper was also financially supported by the EU and the Hungarian government through the project Intensification of the activities of HU-MATHS-IN—Hungarian Service Network of Mathematics for Industry and Innovation under the grant number EFOP-3.6.2-16-2017-00015.

 \bibliographystyle{abbrv}
\bibliography{references}
\appendix
\section{Artificial neural networks}%
\label{sub:artificial_neural_networks}

An artificial neural network (ANN) is a computational model that imitates the workings of a biological neural network. It consists of artificial "neurons" that are loosely inspired by real neurons in the brain. Like other machine learning algorithms, ANNs learn from data to improve their performance.

In particular, these models not only learn how to map input to output, but also how to represent the input itself. This can help explain observed data better.

One type of ANN is called a multilayer perceptron (MLP), or feedforward neural network. This kind of network approximates a function by combining multiple simple functions. MLPs have been widely used in ROM community, for example in the POD-NN method proposed in \cite{ubbiali}. This technique has been used for subsequent investigations in several physical areas, including turbulence \cite{zancanaro}, instability in the convection-dominated problems \cite{hesthaven2019}, aerodynamics \cite{aerostructural-podnn}, and CFD bifurcations \cite{PichiArtificialNeuralNetwork2021}.

 Universal approximation theorem \cite{universal-approximation-theorem} states that for any continuous function $f$, on a compact set $K\subset\mathbb{R}^{n}$, there exists an MLP with one hidden layer that can approximate $f$ uniformly well within any given tolerance $\epsilon>0$. If we use two hidden layers instead of just one, this property extends to all functions - not just continuous ones.

 The network is made up of units called perceptrons, which receive a vector $\boldsymbol{y}_{s}=[y_{s_1},\ldots,y_{s_{m}}]^{\top} \in \mathbb{R}^{m}$ of input values from other neurons and convert it into a single output value $c_{j}$:

\begin{equation}
    c_{j}(\boldsymbol{y}_{s},\boldsymbol{w}_{s,j})= \boldsymbol{w}_{s,j}^{\top} \boldsymbol{y}_{s} ,
\end{equation}

 This output value is then transformed into the excitation state of the neuron through an activation function:

\begin{equation}
    y_{j}=f_{a c t}\left(c_{j};b_{0,j}\right)=f_{act}(\boldsymbol{w}_{s,j}^{\top} \boldsymbol{y}_{j}+b_{0,j}) .
\end{equation}

A typical example of activation function is the \textit{hyperbolic tangent}:
\begin{equation}
f_{a c t}(x)=\frac{e^{x}-e^{-x}}{e^{x}+e^{-x}} .
\end{equation}

 The idea behind MLPs is to connect various perceptrons together in layers, with the first layer representing the input and the last layer representing the output. In between there are hidden layers composed of perceptrons that work according to this model.
The network can be modeled by applying a linear transformation to the vector of the previous layer, followed by a nonlinear transformation represented by the activation function:

\begin{equation}
{\pi}_{NN}\left(\boldsymbol{\mu} ; \Theta\right)\left\{\begin{array}{l}
\boldsymbol{y}_{0}=\boldsymbol{\mu} ; \\
\boldsymbol{y}_{j}=\phi_{j}\left(\mathbf{W}_{j} \boldsymbol{y}_{j-1}+\boldsymbol{b}_{j}\right), \text { for } j=1, \cdots, N  ; \\
\boldsymbol{y}_{out}=\boldsymbol{y}_{N}  .
\end{array}\right.
\end{equation}

Where $N_{j}$ is the number of nodes for the $j-$layer,which is characherized by  the weights $\mathbf{W}_{j} \in \mathbb{R}^{N_{j}\times N_{j-1}}$ and biases $b_{j} \in \mathbb{R}^{N_{j}}$.
We can simplify the neural network by introducing the transformation maps $\pi_{j}(\boldsymbol{y}_{j-1})=\phi_{j}\left(\mathbf{W}_{j} \boldsymbol{y}_{j-1}+\boldsymbol{b_
}{i}\right)$, which allows us to rewrite it as:
\begin{equation}
\label{nn-as-composition}
\pi_{NN}:\left(\boldsymbol{\mu} ; \Theta\right) \mapsto \pi_{N}\left(\cdot ; \mathbf{W}_{N}; \boldsymbol{b}_{N}\right) \circ \ldots \circ \pi_{1}\left(\boldsymbol{\mu} ; \mathbf{W}_{0};\boldsymbol{b}_{1}\right) ,
\end{equation}

where we have collected all the network parameters within $\Theta$.
A supervised learning paradigm is then used to find the optimal parameter vector $\Theta$, i.e., given a set of $N_{\mu}$ input-output pairs $\{\boldsymbol{\mu}_{i},\boldsymbol{y}^{i}\}_{i=1}^{N_{\mu}}$, optimize a loss (or cost) function on the training data:

\begin{equation}
    \label{nn-opt}
    \Theta=\underset{\hat{\Theta}}{\operatorname{argmin}} \ \mathcal{J}(\hat{\Theta}) =\\\underset{\hat{\Theta}}{\operatorname{argmin}} \frac{1}{N_{\mu}} \sum_{i=1}^{N_{\mu}} \mathcal{L}\left(\boldsymbol{y}^{i}, \boldsymbol{y}_{N}^{i} ; \hat{\Theta}\right) ,
\end{equation}
where $\mathcal{J}$ is the loss function and $\mathcal{L}$ the per-sample loss function.

A common choice for measuring the per-sample error is to use the squared euclidean distance, which corresponds to the cost function known as cumulative mean square error (MSE):

\begin{equation}
    \mathcal{J}(\Theta)=\frac{1}{N_{\mu}} \sum_{i=1}^{N_{\mu}} \|\boldsymbol{y}^{i}-\boldsymbol{y}_{N}^{i}\|^{2}_{\mathbb{R}^{N_{r}}} .
\end{equation}

The recent success of neural networks is due in part to the efficient solution of optimization problem \eqref{nn-opt}. The most popular method for doing this is called stochastic gradient descent (SGD), which  employs back-propagation \cite{back-prop} to compute the partial derivatives of the loss function with respect to the parameters. We refer readers to \cite{goodfellow} for an in-depth description of the different techniques that can be used to solve Problem \eqref{nn-opt}.

\section{Source emission database}%
\label{sec:source_emission_database}
The source term, the pollution emission is modelled using the integrated urban air pollution dispersion modelling framework developed by Horv\'ath et al. \cite{Horvath2016410}. The framework attaches the COPERT model to the SUMO traffic simulator software \cite{Sumo}. SUMO (Simulation of Urban Mobility) is an open-source traffic simulator that models road networks, vehicles, and traffic flow behavior in urban areas. It can simulate various traffic scenarios and analyze the environmental impact of transportation. SUMO interfaces with other simulation tools and provides visualizations and statistical data to analyze simulation results. It has been used in transportation planning, traffic management, and the development of autonomous vehicles.
In SUMO, a microscopic traffic model is used, which simulates the behavior of individual vehicles in the traffic network. The physical characteristics of the road, such as speed limits and lane configurations are taken into consideration along with driver behavior like acceleration and braking patterns, lane changes and reaction times. This way realistic traffic flow behavior is obtained making it possible to evaluate the impact of different traffic scenarios on factors such as emissions and fuel consumption.
The COPERT (COmputer Program to calculate Emissions from Road Transport) emission model is a widely-used software tool for estimating emissions from road vehicles \cite{Copert}. It takes various factors  into account that influence vehicle emissions, such as vehicle type, engine size, fuel type, driving cycle, and age. It uses these factors to estimate emissions of pollutants such as carbon monoxide (CO), nitrogen oxides (NOx), particulate matter (PM), and greenhouse gases such as carbon dioxide (CO2).
For this investigation, a small neighborhood of the Bologna University was considered. Using the network tool in SUMO, the road configuration was downloaded and used for simulation. Traffic was generated for every day based on random trip generation. A specific vehicle entered the road network, and took a random exit at every intersection. Only trucks, buses and passenger cars were used. The input of the model is the arrival rate, the average time between two vehicles entering the network. This fluctuated between 80 and 2 seconds, depending on the time of day. Additionally, weekday and monthly factors were considered in a range of 0.6 to 1.0.
The resulting traffic data was aggregated on 5-minute time frames and 20 meter road segments. All vehicles were noted with type and speed. Then, the COPERT model was used to calculate NOx emission in ug/s per road segment and 5-minute frame based on vehicle type and speed. The resulting emission is mapped onto the prepared grid on a per cell basis.

\section{Some notes on ITHACA-FV}
\label{sec:ithaca}
ITHACA-FV is an implementation in OpenFOAM of several reduced order modelling techniques. ITHACA-FV is designed for the ESI version (https://www.openfoam.com/) of OpenFOAM (v1812 onward) but it can be easily adapted also to other versions of OpenFOAM.

ITHACA-FV can also be used as a basis for more advanced projects that would like to assess the capability of reduced order models in their existing OpenFOAM-based software, thanks to the availability of several reduced order methods and algorithms.

Linear and non-linear algebra operations which are not already implemented in OpenFOAM are performed with the external library Eigen. The source code of Eigen 3.4 is provided together with ITHACA-FV as a git submodule and is located in the src/thirdyparty/Eigen folder. For the EigenValue decomposition it is also possible to rely on the Spectra-0.7.0 library and the source code is provided in the src/thirdyparty/spectra folder. Numerical optimization can be performed using the external library OptimLib and the header based source code is provided in the src/thirdyparty/OptimLib folder.

ITHACA-FV has been tested on ubuntu 16.04, CentOS 7, and ArchLinux; it can be easily compiled on any linux distribution with a compiled version of OpenFOAM.

Note that in order to save time and reduce the computational cost, ITHACA-FV skips the offline computations when the offline solutions are already available.
\end{document}